\documentclass[11pt]{article}

\usepackage{amsmath,amssymb,enumitem}
\usepackage{pstricks,pst-node,xy}
\xyoption{all}

\usepackage{xr}
\numberwithin{equation}{section}
\newtheorem{thm}{Theorem}[section] 
\newtheorem{prp}[thm]{Proposition}
\newtheorem{lmm}[thm]{Lemma}   

\newtheorem{dfn}[thm]{Definition}

\def\ov#1{\overline{#1}}
\def\ti#1{\tilde{#1}}
\def\wt#1{\widetilde{#1}}
\def\e_ref#1{(\ref{#1})}
\def\smsize#1{\begin{small}#1\end{small}}

\def\lra{\longrightarrow}
\def\Lra{\Longrightarrow}
\def\Llra{\Longleftrightarrow}

\def\al{\alpha}
\def\be{\beta}
\def\ga{\gamma}
\def\de{\delta}

\def\si{\sigma}
\def\th{\theta}

\def\Ga{\Gamma}
\def\Om{\Omega}
\def\Th{\Theta}

\def\i{\infty}
\def\hb{\hbar}

\def\cA{\mathcal A}

\def\C{\mathbb C}
\def\cC{\mathcal C}
\def\d{\mathfrak d}
\def\D{\mathfrak D}
\def\cD{\mathcal D}
\def\E{\mathbf e}

\def\H{\mathcal H}
\def\I{\mathfrak i}

\def\cM{\mathcal M}
\def\M{\mathfrak M}
\def\N{\mathcal N}
\def\P{\mathbb P}

\def\Pn{\mathbb P^n}

\def\O{\mathcal O}
\def\Q{\mathbb Q}
\def\T{\mathbb T}
\def\X{\mathfrak X}
\def\U{\mathfrak U}
\def\V{\mathcal V}
\def\cY{\mathcal Y}
\def\Z{\mathbb Z}
\def\cZ{\mathcal Z}

\def\Aut{\textnormal{Aut}}
\def\Edg{\textnormal{Edg}}

\def\ev{\textnormal{ev}}

\def\mod{\textnormal{mod~}}
\def\Res{\textnormal{Res}}
\def\Sym{\textnormal{Sym}}

\def\val{\textnormal{val}}
\def\Ver{\textnormal{Ver}}

\begin{document}

\thispagestyle{empty}

\title{Genus-Zero Two-Point Hyperplane Integrals\\
in the Gromov-Witten Theory}
\author{Aleksey Zinger\thanks{Partially supported by a Sloan fellowship 
and DMS Grant 0604874}}
\date{\today}
\maketitle

\begin{abstract}
\noindent
In this paper we compute certain two-point integrals over 
a moduli space of stable maps into projective space.
Computation of one-point analogues of these integrals 
constitutes a proof of mirror symmetry for genus-zero one-point 
Gromov-Witten invariants of projective hypersurfaces. 
The integrals computed in this paper constitute a significant
portion in the proof of mirror symmetry for genus-{\it{one}} 
GW-invariants completed in a separate paper.
These integrals also provide explicit mirror formulas 
for genus-zero {\it two}-point 
GW-invariants of projective hypersurfaces.
The approach described in this paper leads to a reconstruction
algorithm for all genus-zero GW-invariants of
projective hypersurfaces.
\end{abstract}

\tableofcontents

\section{Introduction}
\label{intro_sec}

\subsection{Background and Motivation} 
\label{backgr_subs}

\noindent
The theory of Gromov-Witten invariants has been greatly influenced by 
its interactions with string theory.
In particular, the mirror symmetry principle has led to completely 
unexpected predictions concerning GW-invariants of Calabi-Yau manifolds.
The original prediction of~\cite{CDGP} for the genus-zero GW-invariants
of a quintic threefold was verified about ten years ago in a variety of ways
in~\cite{Ber}, \cite{Ga}, \cite{Gi}, \cite{Le}, and~\cite{LLY}.
The 1993 prediction of~\cite{BCOV} for the genus-one GW-invariants of a quintic 
threefold is verified in~\cite{bcov1}, using the results of this paper.\\

\noindent
The proof of the genus-zero mirror symmetry  for a projective
hypersurface~$X$  essentially consists of 
computing certain equivariant integrals on moduli spaces
$\ov\M_{0,m}(\P^{n-1},d)$ of stable degree-$d$ maps from genus-zero curves
with $m$ marked points into~$\P^{n-1}$.
While the integrals appearing in Chapters~29 and~30 of~\cite{MirSym} are 
over $\ov\M_{0,2}(\P^{n-1},d)$, the integrands involve only one marked point.
For this reason, such integrals can be easily expressed in terms of integrals
on~$\ov\M_{0,1}(\P^{n-1},d)$ and 
determine genus-zero one-point GW-invariants of~$X$; see~\e_ref{Z1ptdfn_e1} below.
In this paper we compute integrals on~$\ov\M_{0,2}(\P^{n-1},d)$ with integrands
involving both marked points.
These integrals in a sense correspond to arithmetic genus one
and indeed constitute a significant portion of the proof of 
mirror symmetry for genus-{\it{one}} GW-invariants in~\cite{bcov1}.
Theorem~\ref{main_thm} also provides closed mirror formulas for
genus-zero {\it two}-point GW-invariants of~$X$, including with descendants.
At the end of Subsection~\ref{outline_subs}, 
we describe the issue arising for integrals with more marked points,
a potential way of addressing it, and a reconstruction algorithm
for genus-zero GW-invariants of $X$ with descendants.\\

\noindent
Let $\U$ be the universal curve over $\ov\M_{0,m}(\P^{n-1},d)$,
with structure map~$\pi$ and evaluation map~$\ev$:
$$\xymatrix{\U \ar[d]^{\pi} \ar[r]^{\ev} & \P^{n-1} \\
\ov\M_{0,m}(\P^{n-1},d).}$$
In other words, the fiber of $\pi$ over $[\cC,f]$ is the curve~$\cC$ 
with $m$ marked points, while 
$$\ev\big([\cC,f;z]\big)=f(z) \qquad\hbox{if}\quad 
z\!\in\!\cC.$$
If $a$ is a positive integer, the sheaf 
$$\pi_*\ev^*\O_{\P^{n-1}}(a)\lra \ov\M_{0,m}(\P^{n-1},d)$$
is locally free. 
We denote the corresponding vector bundle by 
$\V_0$.\footnote{The fiber of $\V_0$ over a point $[\cC,f]\!\in\!\ov\M_{0,m}(\P^{n-1},d)$
is $H^0(\cC;f^*\O_{\P^{n-1}}(a))/\Aut(\cC,f)$.}
Its euler class, $e(\V_0)$, relates genus-zero GW-invariants of 
a degree-$a$ hypersurface in $\P^{n-1}$ to genus-zero GW-invariants of~$\P^{n-1}$;
see Section~26.1 in~\cite{MirSym}.
For each $i\!=\!1,\ldots,m$, there is a well-defined bundle map
$$\wt\ev_i\!:\V_0\lra\ev_i^*\O_{\P^{n-1}}(a), \qquad
\wt\ev_i\big([\cC,f;\xi]\big)=\big[\xi(x_i(\cC))\big],$$
where $x_i(\cC)$ is the $i$th marked of $\cC$.
Since it is surjective, its kernel is again a vector bundle.
Let
$$\V_0'=\ker\wt\ev_1\lra \ov\M_{0,m}(\P^{n-1},d)
\qquad\hbox{and}\qquad
\V_0''=\ker\wt\ev_2\lra \ov\M_{0,m}(\P^{n-1},d),$$
whenever $m\!\ge\!1$ and $m\!\ge\!2$, 
respectively.\footnote{In Chapters 29 and 30 of~\cite{MirSym},
the roles of the marked points $1$ and~$2$ in~\e_ref{Z1ptdfn_e1} are switched;
the analogues of $\V_0$ and $\V_0'$ over $\ov\M_{0,2}(\P^n,d)$ are denoted by
$E_{0,d}$ and $E_{0,d}'$, respectively.}\\

\noindent
The standard action of the $n$-torus $\T$ on $\P^{n-1}$ induces $\T$-actions on 
$\ov\M_{0,m}(\P^{n-1},d)$, $\U$, $\V_0$, $\V_0'$, and~$\V_0''$; 
see Subsections~\ref{equivcoh_subs} and~\ref{local_subs1} for details on
equivariant cohomology.
In particular, $\V_0$, $\V_0'$, and~$\V_0''$ have well-defined equivariant euler classes
$$\E(\V_0),\E(\V_0'),\E(\V_0'')\in H_{\T}^*\big(\ov\M_{0,m}(\P^{n-1},d)\big).$$
These classes are related by
\begin{equation}\label{eulerclass_e}
\E(\V_0)=a\,\ev_1^*(x)\,\E(\V_0')=a\,\ev_2^*(x)\,\E(\V_0''),
\end{equation}
where $x\!\in\!H_{\T}^*(\P^{n-1})$ is the equivariant hyperplane class.
For each $i\!=\!1,2,\ldots,m$, there is also a well-defined equivariant $\psi$-class,
$$\psi_i\in H_{\T}^2\big(\ov\M_{0,m}(\P^{n-1},d)\big),$$
the first chern of the vertical cotangent line bundle of $\U$ 
pull-backed to $\ov\M_{0,m}(\P^{n-1},d)$ by the section
$$\ov\M_{0,m}(\P^{n-1},d)\lra\U, \qquad [\cC,f]\lra\big[\cC,f;x_i(\cC)\big].$$
Since $\ov\M_{0,m}(\P^{n-1},d)$ is a smooth stack (orbifold), there is 
an integration-along-the-fiber homomorphism
$$\int_{\ov\M_{0,m}(\P^{n-1},d)}\!: H_{\T}^*\big(\ov\M_{0,m}(\P^{n-1},d)\big)
\lra H_{\T}^*\approx\Q[\al_1,\ldots,\al_n].$$
For each $i\!=\!1,2,\ldots,n$, let
$$\phi_i\in H_{\T}^*(\P^{n-1})
\approx\Q[x,\al_1,\ldots,\al_n]\big/(x\!-\!\al_1)\ldots(x\!-\!\al_n)$$
be the equivariant Poincare dual of the $i$th fixed point $P_i\!\in\!\P^{n-1}$.
Let
$$\Q_{\al}\equiv\Q(\al_1,\ldots,\al_n)$$
denote the field of fractions in $\al_1,\ldots,\al_n$.
For $a\!=\!1,2,\ldots,n$, an explicit algebraic formula for  
\begin{equation}\label{Z1ptdfn_e1}\begin{split}
\cZ(\hb,\al_i,u) &\equiv  1+\sum_{d=1}^{\i}u^d
\int_{\ov\M_{0,2}(\P^{n-1},d)}\frac{\E(\V_0')}{\hb\!-\!\psi_1}\ev_1^*\phi_i\\
&=\hb^{-1}\Bigg(\hb+ \sum_{d=1}^{\i}u^d
\int_{\ov\M_{0,1}(\P^{n-1},d)}\frac{\E(\V_0')}{\hb\!-\!\psi_1}\ev_1^*\phi_i\Bigg)
\in \Big(\Q_{\al}[[\hb^{-1}]]\Big)\big[\big[u\big]\big]
\end{split}\end{equation}
is confirmed in Chapters~29 and~30 of~\cite{MirSym}.
The equality in~\e_ref{Z1ptdfn_e1} is a straightforward consequence of 
the string relation for GW-invariants; see Section~26.3 in~\cite{MirSym}.\\

\noindent
One of the ingredients in genus-{\it one} localization computations is 
a two-pointed version of~\e_ref{Z1ptdfn_e1}:
\begin{equation}\label{Z2ptdfn_e1}\begin{split}
\wt\cZ(\hb_1,\hb_2,\al_i,\al_j,u) \equiv&  
\frac{a\,\al_i}{\hb_1\!+\!\hb_2}\prod_{k\neq i}(\al_j\!-\!\al_k)\\
&+\sum_{d=1}^{\i}u^d
\int_{\ov\M_{0,2}(\P^{n-1},d)}\frac{\E(\V_0)\ev_1^*\phi_i\ev_2^*\phi_j}
{(\hb_1\!-\!\psi_1)(\hb_2\!-\!\psi_2)}
\in \Big(\Q_{\al}\big[\big[\hb_1^{-1},\hb_2^{-1}\big]\big]\Big)\big[\big[u\big]\big].
\end{split}\end{equation}
Note that the term of degree zero in $u$ above is symmetric in 
$(\hb_1,\al_i)$ and $(\hb_2,\al_j)$, just as are the positive-degree terms.
In turn, $\wt\cZ(\hb_1,\hb_2,\al_i,\al_j,u)$ can be determined from the power series 
\begin{equation}\label{Z2wptdfn_e1}
\cZ_p(\hb,\al_i,u) \equiv \al_i^{p+1}+\sum_{d=1}^{\i}u^d
\int_{\ov\M_{0,2}(\P^{n-1},d)}\frac{\E(\V_0'')\ev_2^*x^{p+1}}{\hb\!-\!\psi_1}
\ev_1^*\phi_i \in \Big(\Q_{\al}[[\hb^{-1}]]\Big)\big[\big[u\big]\big]
\end{equation}
with $p\!=\!-1,0,\ldots,n\!-\!1$.
The seemingly unfortunate choice of indexing is partly motivated by the central role played
by the power series $\cZ(\hb,\al_i,u)$ defined in~\e_ref{Z1ptdfn_e1} and the simple relation
$$\cZ_0(\hb,\al_i,u)=\al_i\cZ(\hb,\al_i,u),$$
which follows from~\e_ref{eulerclass_e}, along with 
\e_ref{ABothm_e}, \e_ref{phidfn_e}, and~\e_ref{restrmap_e}.
As shown in this paper,
\begin{equation}\label{Z2ptform_e}
\wt\cZ(\hb_1,\hb_2,\al_i,\al_j,u)
=\frac{a}{\hb_1\!+\!\hb_2}\sum_{p+q+r=n-1}\!\!\!\!\!\!
(-1)^r\si_r \cZ_p(\hb_1,\al_i,u)\cZ_{q-1}(\hb_2,\al_j,u),
\end{equation}
where $\si_p$ is the $p$th elementary symmetric polynomial in $\al_1,\ldots,\al_n$;
see Theorem~\ref{main_thm}.\\

\noindent
{\it Remark 1:} The right-hand side of~\e_ref{Z2ptform_e} is in fact symmetric in
$(\hb_1,\al_i)$ and $(\hb_2,\al_j)$, because
\begin{equation}\label{symmrel_e}
\cZ_{n-1}(\hb,\al_i,u)-\si_1\cZ_{n-2}(\hb,\al_i,u)+\ldots
+(-1)^n\si_n\cZ_{-1}(\hb,\al_i,u)=0.
\end{equation}
The reason for this relation is explained in Subsection~\ref{outline_subs}.\\

\noindent
{\it Remark 2:} We will see in Subsection~\ref{mainthm_subs} that 
the power series $\cZ_p(\hb,\al_i,u)$ can be represented by elements
of $\Q_{\al}(\hb)[[u]]$.
The relation \e_ref{Z2ptform_e} might then suggest  that the corresponding
element of $\Q_{\al}(\hb_1,\hb_2)[[u]]$ representing $\wt\cZ(\hb_1,\hb_2,\al_i,\al_j,u)$
has a simple pole at $\hb_1\!=\!-\hb_2$.
In fact, there is no pole at $\hb_1\!=\!-\hb_2$, except in degree zero.
This is immediate from the localization formula~\e_ref{ABothm_e}; 
see also Subsection~\ref{local_subs1}.\\

\noindent
The power series~\e_ref{Z2wptdfn_e1} encode genus-zero two-point GW-invariants
of a degree-$a$ hypersurface in~$\P^{n-1}$ with constraints coming from~$\P^{n-1}$.
Thus, Theorem~\ref{main_thm} provides mirror formulas for such invariants;
the coefficients $\ti{C}_{p,q}^{(r)}$ are ``purely equivariant'' 
and are irrelevant for this purpose.
In the table below, we give the first ten genus-zero two-point BPS numbers,
defined from GW-invariants by equation~(2) in~\cite{KP}, 
for the degree-$6$ hypersurface in~$\P^7$.
These numbers are integers as predicted by Conjecture~1 in~\cite{KP}.
In fact, we have used the first statement of Theorem~\ref{main_thm},
along with a computer program, to confirm this conjecture for 
all degree-$d$ two-point BPS counts in a degree-$n$ hypersurface 
$X_n$ in~$\P^{n-1}$ with $n\!\le\!10$ and $d\!\le\!20$.\\

\begin{tabular}{l|l}
\hline
degree $d$& BPS curve count through 2 codim-2 linear subspaces in $X_7$\\
\hline
1& 1707797\\
2& 510787745643\\
3& 222548537108926490\\
4& 113635631482486991647224\\
5& 63340724462384110502639024265\\
6& 37325795060717360046547665187418254\\
7& 22857028298936684292245509537579343818647\\
8& 14395953469762596243721601709186933042635134584\\
9& 9263611884884554518268724722981763557936573405648178\\
10& 6062677702410680024315392235188823274104219383883410807999\\
\hline
\end{tabular}\\

\vspace{.1in}

\noindent
The explicit expressions of Section 29.1 in \cite{MirSym} for the power 
series $\cZ(\hb,\al_i,u)$ have very different forms for $a\!<\!n$ and $a\!=\!n$.
The $a\!=\!n$ case is the most interesting one and corresponds to Calabi-Yau hypersurfaces.
As the power series $\cZ(\hb,\al_i,u)$ are central to our computation of 
$\cZ_p(\hb,\al,u)$ and $\wt\cZ(\hb_1,\hb_2,\al_i,\al_j,u)$, 
for the purposes of the explicit expressions preceding Theorem~\ref{main_thm}
in the next subsection we consider only the case $a\!=\!n$.\footnote{In other words,
one may to choose set $a\!=\!n$ for the rest of the paper.
However, the statements of Lemmas~\ref{recgen_lmm} and~\ref{polgen_lmm} and 
their proofs are valid for all~$a$.
Therefore, the proofs of~\e_ref{Z2ptform_e} and \e_ref{symmrel_e} are valid as well.
The same is the case with Theorem~\ref{main_thm} if the power series $\cY_{-1}$ and
$\cY$ are chosen appropriately.} 
This is also the case used in~\cite{bcov1}.

\subsection{Main Theorem}
\label{mainthm_subs}

\noindent
The essence of mirror symmetric predictions for Gromov-Witten invariants
is that these invariants (and relatedly $\cZ(\hb,\al_i,u)$)
can be expressed in terms of certain hypergeometric series.
In this subsection, we define these series and then express
$\cZ(\hb,\al_i,u)$, $\cZ_p(\hb,\al_i,u)$, and $\wt\cZ(\hb_1,\hb_2,\al_i,\al_j,u)$
in terms of them.\\

\noindent
Let $n$ be a positive integer.
For each $q\!=\!0,1,\ldots$, define $I_{0,q}(t)$ by 
\begin{equation}\label{Ifuncdfn_e}
\sum_{q=0}^{\i} I_{0,q}(t)w^q \equiv e^{wt}\sum_{d=0}^{\i}e^{dt}
\frac{\prod_{r=1}^{r=nd}(nw\!+\!r)}{\prod_{r=1}^{r=d}(w\!+\!r)^n}.
\end{equation}
Each $I_{0,q}(t)$ is a degree-$q$ polynomial in $t$ 
with coefficients that are power series in~$e^t$.
For example,
\begin{equation}\label{HGexamp_e}
I_0(t)=1+\sum_{d=1}^{\i}e^{dt}\frac{(nd)!}{(d!)^n}
\qquad\hbox{and}\qquad
I_1(t)=tI_0(t)+\sum_{d=1}^{\i}e^{dt}
\bigg(\frac{(nd)!}{(d!)^n}\sum_{r=d+1}^{nd}\!\frac{n}{r}\bigg).
\end{equation}\\

\noindent
For $p,q\!\in\!\Z^+$ with $q\!\ge\!p$, let
\begin{equation}\label{Tden_e}
I_{p,q}(t)=\frac{d}{dt}\bigg(\frac{I_{p-1,q}(t)}{I_{p-1,p-1}(t)}\bigg).  
\end{equation}
It is straightforward to check that each of the ``diagonal'' terms $I_{p,p}(t)$ 
is a power series in $e^t$ with constant term~$1$, whenever it is defined;
see \cite{ZaZ}, for example. 
Thus, the division in~\e_ref{Tden_e} is well-defined for all~$p$.
Let
\begin{equation}\label{mirmap_e}
T=\frac{I_{0,1}(t)}{I_{0,0}(t)}.
\end{equation}
By~\e_ref{HGexamp_e}, the map $t\!\lra\!T$ is a change of variables;
it will be called the {\tt mirror map}.
If $p\!\in\!\bar\Z^+$ and $\cY(\hbar,x,e^t)$ is a  power series in $e^t$ with coefficients 
that are functions of a complex variable $\hbar$ and possibly some other variable $x$, let
\begin{equation}\label{Tderivdfn_e}
\D^p\cY(\hbar,x,t)=e^{-xt/\hb}
\bigg\{\frac{\hb}{I_{p,p}(t)}\frac{d}{dt}\Big\}\ldots
\Big\{\frac{\hb}{I_{1,1}(t)}\frac{d}{dt}\Big\}
\big(e^{xt/\hb}\cY(\hb,x,e^t)\big).
\end{equation}\\

\noindent
We define
$$\cY(\hbar,x,e^t)=I_{0,0}(t)^{-1}x\sum_{d=0}^{\i}e^{dt}
\frac{\prod_{r=1}^{r=nd}(nx\!+\!r\hb)}
{\prod_{r=1}^{r=d}\prod_{k=1}^{k=n}(x\!-\!\al_k\!+\!r\hb)}
\in\Q_{\al}(\hb,x)\big[\big[e^t\big]\big]\Big/\prod_{k=1}^{k=n}(x\!-\!\al_k).$$
Expanding $e^{xt/\hb}\cY(\hbar,x,e^t)$ as a power series in~$\hb^{-1}$,
we obtain
\begin{equation}\label{Yexpand_e}
e^{xt/\hb}\cY(\hbar,x,e^t)=x
\sum_{q=0}^{\i}\bigg(\sum_{r=0}^{r=q}C_{0,q}^{(r)}(t)x^{q-r}\bigg)\hb^{-q},
\end{equation}
where $C_{0,q}^{(r)}(t)$ is a degree-$r$ symmetric polynomial in $\al_1,\ldots,\al_n$
with coefficients in $\Q[t][[e^t]]$.
For example, 
\begin{equation}\label{Yexpand_e2}
C_{0,q}^{(0)}(t)=I_{0,q}(t)\big/I_{0,0}(t) \qquad\hbox{and}\qquad
C_{0,1}^{(1)}(t)=\si_1I_{0,0}(t)^{-1}\sum_{d=1}^{\i}e^{dt}
\bigg(\frac{(nd)!}{(d!)^n}\sum_{r=1}^{n}\!\frac{1}{r}\bigg).
\end{equation}
The main conclusion of Section 30.4 in~\cite{MirSym} is that the power series 
$\cZ(\hb,\al_i,e^T)$ defined in~\e_ref{Z1ptdfn_e1} is the evaluation~of
\begin{equation}\label{Z1ms_e}
\cZ(\hb,x,e^T)= 
e^{(t-T)x/\hb}e^{-C_{0,1}^{(1)}(t)/\hb}\cY(\hbar,x,e^t)\in
\Q_{\al}(\hb,x)\big[\big[e^t\big]\big]\Big/\prod_{k=1}^{k=n}(x\!-\!\al_k)
\end{equation}
at $x\!=\!\al_i$, if $T$ and $t$ are related by the mirror map~\e_ref{mirmap_e}.\\

\noindent
The power series $\cZ_p(\hb,\al_i,u)$ and $\wt\cZ(\hb_1,\hb_2,\al_i,\al_j,u)$ defined
in~\e_ref{Z2wptdfn_e1} and~\e_ref{Z2ptdfn_e1}, respectively,
are also evaluations of certain power series
\begin{equation*}\begin{split}
&\cZ_p(\hb,x,u)\in\Q_{\al}(\hb,x)\big[\big[u\big]\big]
\Big/\prod_{k=1}^{k=n}(x\!-\!\al_k)  \qquad\hbox{and}\\
&\wt\cZ(\hb_1,\hb_2,x_1,x_2,u)\in
\Q_{\al}(\hb_1,\hb_2,x_1,x_2)\big[\big[u\big]\big]\Big/
\prod_{k=1}^{k=n}(x_1\!-\!\al_k)(x_2\!-\!\al_k)
\end{split}\end{equation*}
that have a mirror transform shape analogous to~\e_ref{Z1ms_e}.
Let
\begin{equation}\label{Ym1dfn_e}
\D^{-1}\cY(\hbar,x,e^t)\equiv \sum_{d=0}^{\i}e^{dt}
\frac{\prod_{r=0}^{nd-1}(nx\!+\!r\hb)}{
\prod_{r=1}^{r=d}\prod_{k=1}^{k=n}(x\!-\!\al_k\!+\!r\hb)}
\in\Q_{\al}(\hb,x)\big[\big[e^t\big]\big]\Big/\prod_{k=1}^{k=n}(x\!-\!\al_k).
\end{equation}

\begin{thm}
\label{main_thm}
There exist  $\ti{C}_{p,q}^{(r)}\!\in\!\Q_{\al}\big[[e^t]\big]$, 
with $p\!\ge\!r\!\ge\!1$ and $p\!-\!r\!\ge\!q\!\ge\!0$,
such that the coefficients of the powers of $e^t$ in $\ti{C}_{p,q}^{(r)}$ 
are degree-$r$ symmetric polynomials and the power series defined in~\e_ref{Z2wptdfn_e1}
are given~by
\begin{equation}\label{main_e1}\begin{split}
\cZ_p(\hb,x,e^T)&= 
e^{(t-T)x/\hb}e^{-C_{0,1}^{(1)}(t)/\hb}\cY_p(\hbar,x,e^t),
\qquad\hbox{where}\\
\cY_p(\hb,x,e^t)&\equiv\D^p\cY(\hbar,x,e^t)
+\sum_{r=1}^{r=p}\sum_{q=0}^{p-r}\ti{C}_{p,q}^{(r)}(e^t)\hb^{p-r-q}\D^q\cY(\hbar,x,e^t),
\end{split}\end{equation}
if $T$ and $t$ are related by the mirror map~\e_ref{mirmap_e}.
Furthermore, the power series defined in~\e_ref{Z2ptdfn_e1} are given~by
\begin{equation}\label{main_e2}
\wt\cZ(\hb_1,\hb_2,x_1,x_2,u)
=\frac{a}{\hb_1\!+\!\hb_2}\sum_{p+q+r=n-1}\!\!\!\!\!\!
(-1)^r\si_r \cZ_p(\hb_1,x_1,u)\cZ_{q-1}(\hb_2,x_2,u).
\end{equation}\\
\end{thm}
 
\noindent
{\it Remark:} We note that by~\e_ref{main_e1} and~\e_ref{main_e2}
\begin{equation*}\begin{split}
\wt\cZ(\hb_1,\hb_2,x_1,x_2,e^T)&= 
e^{(t-T)(x_1/\hb_1+x_2/\hb_2)}e^{-C_{0,1}^{(1)}(t)(\hb_1^{-1}+\hb_2^{-1})}
\wt\cY(\hbar_1,\hb_2,x_1,x_2,e^t),
\qquad\hbox{where}\\
\wt\cY(\hb_1,\hb_2,x_1,x_2,u)&\equiv 
\frac{a}{\hb_1\!+\!\hb_2}\sum_{p+q+r=n-1}\!\!\!\!\!\!
(-1)^r\si_r \cY_p(\hb_1,x_1,u)\cY_{q-1}(\hb_2,x_2,u). 
\end{split}\end{equation*}
In other words, $\wt\cZ$ is the same transform of $\wt\cY$ in both $(x_1,\hb_1)$ and 
$(\hb_2,x_2)$  as $\cZ$ and $\cZ_p$ are of $\cY$ and $\cY_p$ in~$(\hb,x)$.\\

\noindent
The only relevant property of the power series $\ti{C}_{p,q}^{(r)}$
for the purposes of the genus-one localization computations in~\cite{bcov1} is that 
the $e^t$-coefficients of $\ti{C}_{p,q}^{(r)}$ lie in the ideal generated 
by $\si_1,\ldots,\si_{n-1}$ if $p\!\le\!n\!-\!1$.
This is automatic in the case of Theorem~\ref{main_thm}, since each of these
coefficients is a symmetric polynomial in $\al_1,\ldots,\al_n$ of a positive 
degree $r\!\le\!n\!-\!1$.
The approach of~\cite{bcov1} suggests that these coefficients are likely to be
irrelevant in many other localization computations as well.
Nevertheless, they are described inductively in this paper in the process
of proving the first statement of Theorem~\ref{main_thm}.

\subsection{Outline of the Proof} 
\label{outline_subs}

\noindent
The proof of~\e_ref{Z1ms_e} in Chapter 30 of \cite{MirSym} essentially
consists of showing that
\begin{enumerate}[label=(S\arabic*)]
\item $\cY(\hb,x,u)$ and $\cZ(\hb,x,u)$ satisfy a certain recursion on the $u$-degree;
\item $\cY(\hb,x,u)$ and $\cZ(\hb,x,u)$ satisfy a certain self-polynomiality condition (SPC);
\item\label{equalmodS_item}
the two sides of~\e_ref{Z1ms_e}, viewed as a powers series in $\hb^{-1}$, 
agree mod~$\hb^{-2}$; 
\item if $Y(\hb,x,u)$ satisfies the recursion and the SPC,
so do certain transforms of~$Y(\hb,x,u)$;
\item\label{uniqS_item}
if $Y(\hb,x,u)$ satisfies the recursion and the SPC, it is determined by its 
``mod $\hb^{-2}$ part''.
\end{enumerate}
For the purposes of these statements, in particular \ref{equalmodS_item} 
and~\ref{uniqS_item}, we assume that
$$\cY(\hb,x,u),\cZ(\hb,x,u),Y(\hb,x,u)
\in\Q_{\al}(\hb,x)\big[\big[u\big]\big]\Big/\prod_{k=1}^{k=n}(x\!-\!\al_k).$$
For example, \ref{uniqS_item} means 
\begin{equation*}\begin{split}
&Y(\hb,\al_i,u)\cong\bar{Y}(\hb,\al_i,u) ~~~(\mod \hb^{-2})
\quad\forall\,i\!=\!1,2,\ldots,n\\
&\hspace{2in}\Lra\qquad
Y(\hb,\al_i,u)=\bar{Y}(\hb,\al_i,u) \quad\forall\,i\!=\!1,2,\ldots,n.
\end{split}\end{equation*}\\

\noindent
The proof of~\e_ref{main_e1} in this paper essentially consists of showing  that
\begin{enumerate}[label=(M\arabic*)]
\item\label{recM_item}
 $\cY_p(\hb,x,u)$ and $\cZ_p(\hb,x,u)$ satisfy the recursion~\e_ref{recur_dfn_e};
\item\label{polM_item}
$\big(\cY(\hb,x,u),\cY_p(\hb,x,u)\big)$ and $\big(\cZ(\hb,x,u),\cZ_p(\hb,x,u)\big)$ 
satisfy the mutual polynomiality condition (MPC) of Lemma~\ref{Phistr_lmm2};
\item\label{equalmodM_item} the two sides of~\e_ref{main_e1}, 
 viewed as a powers series in $\hb^{-1}$,  agree mod~$\hb^{-1}$; 
\item\label{transM_item} 
if $Z(\hb,x,u)$ satisfies~\e_ref{recur_dfn_e} and the MPC with respect to $Y(\hb,x,u)$,
the transforms of $Z(\hb,x,u)$ of Lemma~\ref{Phistr_lmm4} satisfy~\e_ref{recur_dfn_e} 
and the MPC with respect to appropriate transforms of $Y(\hb,x,u)$;
\item\label{uniqM_item} if the $u$-constant term of $Y(\hb,x,u)$ is independent of $\hb$ and
nonzero and $Z(\hb,x,u)$ satisfies~\e_ref{recur_dfn_e} and the MPC with respect to $Y(\hb,x,u)$, 
it is determined by its  ``mod $\hb^{-1}$ part''; see Proposition~\ref{uniqueness_prp}.
\end{enumerate}
The statements \ref{equalmodM_item} and \ref{uniqM_item} should be interpreted 
analogously to \ref{equalmodS_item} and~\ref{uniqS_item}.
In other words, the equalities are modulo $\prod_{k=1}^{k=n}(x\!-\!\al_k)$,
or equivalently after the evaluation at each $x\!=\!\al_i$.
Similarly, the requirement on the degree-zero term in $Y(\hb,x,u)$ in~\ref{uniqM_item}
means that it is nonzero even after the evaluation at each $x\!=\!\al_i$.\\

\noindent
The claims of \ref{recM_item} and~\ref{polM_item} concerning $\cZ_p(\hb,x,u)$ are 
special cases of Lemmas~\ref{recgen_lmm} and~\ref{polgen_lmm} below, since
\begin{equation}\label{strrel_e}\begin{split}
\int_{\ov\M_{0,2}(\P^{n-1},d)}\!\!\!
\frac{\E(\V_0'')\ev_2^*\eta}{\hb\!-\!\psi_1}\ev_1^*\phi_i
=\hb\!\!\int_{\ov\M_{0,3}(\P^{n-1},d)}\!\!\!
\frac{\E(\V_0'')\ev_2^*\eta}{\hb\!-\!\psi_1}\ev_1^*\phi_i
\quad\forall\,\eta\!\in\!H_{\T}^*(\P^{n-1}),\,d\!\in\!\Z^+,
\end{split}\end{equation}
by the string relation;  see Section~26.3 in~\cite{MirSym}.

\begin{lmm}
\label{recgen_lmm} For all $m\!\ge\!3$, $\eta_j\!\in\!H_{\T}^*(\P^{n-1})$,
and $\be_j\!\in\!\bar\Z^+$, the power series $\cZ_{\eta,\be}(\hb,x,u)$ defined by
\begin{equation}\label{recgen_e}
\cZ_{\eta,\be}(\hb,x,u)\equiv\sum_{d=0}^{\i}u^d\Bigg(
\int_{\ov\M_{0,m}(\P^{n-1},d)}\!\!\!
\frac{\E(\V_0'')\ev_1^*\phi_i}{\hb\!-\!\psi_1}\prod_{j=2}^{j=m}
\big(\psi_j^{\be_j}\ev_j^*\eta_j\big)\Bigg)
\end{equation}satisfies the recursion~\e_ref{recur_dfn_e}.
\end{lmm}

\begin{lmm}
\label{polgen_lmm} For all $m\!\ge\!3$, $\eta_j\!\in\!H_{\T}^*(\P^{n-1})$,
and $\be_j\!\in\!\bar\Z^+$, the power series $\hb^{m-2}\cZ_{\eta,\be}(\hb,x,u)$,
with $\cZ_{\eta,\be}(\hb,x,u)$ as in Lemma~\ref{recgen_lmm},
satisfies the polynomial condition of Lemma~\ref{Phistr_lmm1} with respect to $\cZ(\hb,x,u)$.
\end{lmm}

\noindent
Our proof of Lemma~\ref{recgen_lmm} is practically identical to the proof in Section~30.1 
of \cite{MirSym} that $\cZ(\hb,x,u)$ satisfies a certain recursion on the
$u$-degree.\footnote{However, the coefficients $C_i^j(d)$ in \e_ref{recur_dfn_e} 
are ``shifts by one'' of the coefficients in the recursion for~$\cZ(\hb,x,u)$.}
The proof of Lemma~\ref{polgen_lmm} is similar to the proof in Section~30.2 of~\cite{MirSym} 
that $\cZ(\hb,x,u)$ satisfies the~SPC.
However, there are differences in how the key idea for the setup used in~\cite{MirSym}
is utilized. 
An explanation of the modifications and a complete justification of their appropriateness
are not very simple to state.
In order to avoid any confusion, we thus give a full account in 
Subsections~\ref{local_subs2} and~\ref{PhiZstrpf_subs}.
As it requires most of what constitutes the proof of the recursivity 
relation~\e_ref{recur_dfn_e}, 
we give a proof of the latter in Subsection~\ref{local_subs1}.\\

\noindent
We are now able to justify~\e_ref{symmrel_e}.
By~\e_ref{Z2wptdfn_e1},
\begin{equation}\label{Z2wmod_e}
\cZ_p(\hb,x,u) \cong x^{p+1} \qquad (\mod \hb^{-1}\big).\footnotemark
\end{equation}
\footnotetext{Such identities will be taken to mean that the two sides are equal
if $x$ is replaced with $\al_i$ for every $i\!=\!1,2,\ldots,n$.}
Along with~\ref{uniqM_item}, these three properties of $\cZ_p(\hb,x,u)$ imply~\e_ref{symmrel_e},
since
\begin{equation*}\begin{split}
\sum_{p+r=n}(-1)^r\si_r\cZ_{p-1}(\hb,x,u)
&\equiv \sum_{p+r=n}(-1)^r\si_rx^p \qquad (\mod \hb^{-1}\big)\\
&=\prod_{k=1}^{k=n}(x\!-\!\al_k).
\end{split}\end{equation*}
The last expression above vanishes at $x\!=\!\al_i$ for all $i\!=\!1,2,\ldots,n$.\\

\noindent
We will check by a direct computation that 
$$\cY_{-1}(\hb,x,u)\equiv\D^{-1}\cY(\hb,x,u)$$ 
satisfies~\e_ref{recur_dfn_e} and the MPC with respect to~$\cY(\hb,x,u)$;
see Subsection~\ref{algcomp_subs}.
By~\e_ref{Ym1dfn_e},
$$\cY_{-1}(\hb,x,u)\cong 1 \qquad (\mod \hb^{-1}\big).$$
Thus, the $p\!=\!-1$ case of~\e_ref{main_e1} follows from \e_ref{Z2wmod_e},
\ref{transM_item}, and~\ref{uniqM_item}.\\

\noindent
We could also verify directly that $\cY_0(\hb,x,u)$ satisfies \e_ref{recur_dfn_e}
and the MPC with respect to~$\cY(\hb,x,u)$. 
Fortunately, this is an immediate consequence of parts~\ref{deriv_ch}
and~\ref{mult_ch} of Lemma~\ref{Phistr_lmm4}, since 
$$\cY_0(\hb,x,e^t)=\frac{1}{I_{0,0}(t)}\bigg\{x+\hb\frac{d}{dt}\bigg\}
\cY_{-1}(\hb,x,e^t).$$
Thus, $\cY_0(\hb,x,u)$ satisfies the two properties because $\cY_{-1}(\hb,x,u)$ does.
On the other hand, by~\e_ref{Yexpand_e} and~\e_ref{Yexpand_e2},
$$\cY_0(\hb,x,u)=x\cY(\hb,x,u)\cong x \qquad (\mod \hb^{-1}\big).$$
Thus, the $p\!=\!0$ case of~\e_ref{main_e1} also follows from \e_ref{Z2wmod_e},
\ref{transM_item}, and~\ref{uniqM_item}.\\

\noindent
The differentiation transform~\ref{deriv_ch} of Lemma~\ref{Phistr_lmm4}
is the only one of the five ``admissible'' transforms that has no analogue
in Chapter~30 of \cite{MirSym}.\footnote{The multiplication by $\hb$, 
i.e.~transform~\ref{hbmult_ch}, is not explicitly mentioned in Chapter~30 of \cite{MirSym},
but its admissibility is nearly immediate from the relevant definitions.}
The admissibility of this transform, along with that of~\ref{mult_ch} of Lemma~\ref{Phistr_lmm4},
implies that $\cY_p(\hb,x,u)$ defined by the second equation in~\e_ref{main_e1}  
satisfies the recursion~\e_ref{recur_dfn_e}
and the MPC of Lemma~\ref{Phistr_lmm2} with respect to $\cY(\hb,x,u)$,
no matter what the coefficients $\ti{C}_{p,q}^{(r)}(u)$ are.
In light of~\e_ref{Z2wmod_e}, the $p\!\ge\!1$ cases of the first equation 
in~\e_ref{main_e1} thus reduce to showing there exist $\ti{C}_{p,q}^{(r)}(u)$ such that
\begin{equation}\label{Ypcond_e}
\cY_p(\hb,x,u)\cong x^{p+1}  \qquad (\mod \hb^{-1}\big).
\end{equation}
This is proved by induction, using~\e_ref{Yexpand_e} and~\e_ref{Yexpand_e2};
see Subsection~\ref{algcomp_subs}.\\

\noindent
The proof of \e_ref{main_e2} follows the same principle.
By the string relation and~\e_ref{eulerclass_e},
$$\wt\cZ(\hb_1,\hb_2,\al_i,\al_j,x_2,u)
=\frac{\hb_1\hb_2}{\hb_1\!+\!\hb_2}\sum_{d=0}^{\i}
\int_{\ov\M_{0,3}(\P^{n-1},d)}\!\!\!
\frac{\E(\V_0'')}{(\hb_1\!-\!\psi_1)(\hb_2\!-\!\psi_2)}\ev_1^*\phi_i\ev_2^*(ax\,\phi_j).$$
Thus, by Lemmas~\ref{recgen_lmm} and~\ref{polgen_lmm},
$(\hb_1\!+\!\hb_2)\wt\cZ(\hb_1,\hb_2,x_1,x_2,u)$ 
satisfies the recursion~\e_ref{recur_dfn_e} and the MPC of Lemma~\ref{Phistr_lmm2} 
with respect to $\cZ(\hb,x,u)$ for $(\hb,x)\!=\!(\hb_1,x_1)$ and $x_2\!=\!\al_j$ fixed.
By symmetry, it also satisfies the two properties for $(\hb,x)\!=\!(\hb_2,x_2)$ and 
$x_1\!=\!\al_i$ fixed.
It is then sufficient to compare the two sides of~\e_ref{main_e2} multiplied 
by $\hb_1\!+\!\hb_2$ modulo~$\hb_1^{-1}$:
\begin{alignat*}{1}
&\big(\hb_1\!+\!\hb_2\big)\wt\cZ(\hb_1,\hb_2,\al_i,\al_j,u) \cong 
a\al_i\!\!\!\!\!\!\!  \sum_{p+q+r=n-1}\!\!\!\!\!\!\!\! (-1)^r\si_r \al_i^p\al_j^q
+\sum_{d=1}^{\i}u^d\int_{\ov\M_{0,2}(\P^{n-1},d)}\!\!\!
\frac{\E(\V_0)\ev_1^*\phi_i\ev_2^*\phi_j}{\hb_2\!-\!\psi_2};\\
&a\!\!\!\!\!\!\!\sum_{p+q+r=n-1}\!\!\!\!\!\!\!\!
(-1)^r\si_r \cZ_p(\hb_1,\al_i,u)\cZ_{q-1}(\hb_2,\al_j,u) \cong
a\!\!\!\!\!\!\!\!\!   \sum_{p+q+r=n-1}\!\!\!\!\!\!\!
(-1)^r\si_r\al_i^{p+1}\cZ_{q-1}(\hb_2,\al_j,u).
\end{alignat*}
In order to see that the two resulting expressions are equal, we compare them modulo~$\hb_2^{-1}$:
\begin{alignat*}{2}
&a\!\!\!\!\!\!\!\sum_{p+q+r=n-1}\!\!\!\!\!\!\!\! (-1)^r\si_r \al_i^{p+1}\al_j^q
+\sum_{d=1}^{\i}u^d\int_{\ov\M_{0,2}(\P^{n-1},d)}\!\!\!
\frac{\E(\V_0)\ev_1^*\phi_i\ev_2^*\phi_j}{\hb_2\!-\!\psi_2}
\cong a\!\!\!\!\!\!\!\sum_{p+q+r=n-1}\!\!\!\!\!\!(-1)^r\si_r \al_i^{p+1}\al_j^q;\\
&a\!\!\!\!\!\!\!\sum_{p+q+r=n-1}\!\!\!\!\!\!
(-1)^r\si_r\al_i^{p+1}\cZ_{q-1}(\hb_2,\al_j,u) \cong
a\!\!\!\!\!\!\!\sum_{p+q+r=n-1}\!\!\!\!\!\!(-1)^r\si_r\al_i^{p+1}\al_j^q.
\end{alignat*}
From this we conclude that the two sides of~\e_ref{main_e2} multiplied by $\hb_1\!+\!\hb_2$
modulo~$\hb_1^{-1}$ are equal. Therefore, the two sides~\e_ref{main_e2} are equal
by~\ref{uniqM_item}.\\

\noindent
Central to this paper are the use of the transforms $\D^p$ in conjunction
with part~\ref{deriv_ch} of Lemma~\ref{Phistr_lmm4} and a desymmetrization of
the approach of Chapter~30 of~\cite{MirSym} to obtain an explicit closed formula
for the integrals in~\e_ref{Z2ptdfn_e1}.
The transforms $\D^p$, combined with the transforms~\ref{mult_ch} and~\ref{hbmult_ch}
of Lemma~\ref{Phistr_lmm4}, make it possible to construct a power series $\wt\cY(\hb,x,u)$, satisfying the recursion~\e_ref{recur_dfn_e} and the MPC with respect to~$\cY(\hb,x,u)$,
that agrees with a pre-specified $\al$-symmetric element of $\Q_{\al}[\hb]\big[\big[u\big]\big]$
modulo~$\hb^{-1}$.
On the other hand, for the purposes of \ref{uniqM_item}, 
it is sufficient to assume that the coefficient of each power of $u$ in $Y$ and $Z$
is a sum of a power series in $\hb^{-1}$ and a polynomial 
in~$\hb$.\footnote{The same is the case with~\ref{uniqS_item}. In fact, 
the condition on the middle term of the recursion used Chapter 30 of \cite{MirSym} 
can also be relaxed as in Definition~\ref{recur_dfn}.}
Using~\ref{uniqM_item} and the last two parts of Lemma~\ref{Phistr_lmm4}, 
a variety of integrals on~$\ov\M_{0,m}(\Pn,d)$ involving $\E(\V_0)$ and 
products of $1/(\hb_j\!-\!\psi_j)$ can be reduced to integrals involving~$\E(\V_0)$ 
and powers $\psi$-classes, with each exponent bounded by~$m\!-\!3$.\\

\noindent
For example, suppose one would like to compute 
\begin{equation}\label{Z3ptdfn_e1}
\wt\cZ^{(3)}(\hb_1,\hb_2,\hb_3,\al_{i_1},\al_{i_2},\al_{i_3},u)
\equiv \sum_{d=0}^{\i}u^d
\int_{\ov\M_{0,3}(\P^{n-1},d)}\frac{\E(\V_0)\ev_1^*\phi_{i_1}\ev_2^*\phi_{i_2}\ev_3^*\phi_{i_3}}
{(\hb_1\!-\!\psi_1)(\hb_2\!-\!\psi_2)(\hb_3\!-\!\psi_3)};
\end{equation}
these integrals may be useful for localization computations in (arithmetic) genus two.
By Lemmas~\ref{recgen_lmm} and~\ref{polgen_lmm} and part~\ref{hbmult_ch} of 
Lemma~\ref{Phistr_lmm4}, $\hb_1\hb_2\hb_3\wt\cZ^{(3)}$ satisfies 
the recursion and the MPC with respect to $\cZ$ for each $(\hb,x)\!=\!(\hb_s,x_s)$.
Thus,  $\wt\cZ^{(3)}$ can be reconstructed\footnote{From \e_ref{Z3ptdfn_e2} and
$\phi_{i_3}\!=\!\prod_{k\neq i_3}(x\!-\!\al_{i_3})$, one obtains
$$\hb_1\hb_2\hb_3\wt\cZ^{(3)}(\hb_1,\hb_2,\hb_3,x_1,x_2,x_3,u)=
\hb_1\hb_2\sum_{p+q+r=n-1}\!\!\!(-1)^r\si_r\wt\cZ_p(\hb_1,\hb_2,x_1,x_2,u)
\cZ_q(\hb_3,x_3,u).$$}
from the ``mod $\hb_3^{-1}$'' part of $\hb_1\hb_2\hb_3\wt\cZ^{(3)}$,
\begin{equation}\label{Z3ptdfn_e2}
\hb_1\hb_2\hb_3\wt\cZ^{(3)}(\hb_1,\hb_2,\hb_3,\al_{i_1},\al_{i_2},\al_{i_3},u)
\cong \hb_1\hb_2\sum_{d=0}^{\i}u^d\int_{\ov\M_{0,3}(\P^{n-1},d)}\!\!\!
\frac{\E(\V_0)\ev_1^*\phi_{i_1}\ev_2^*\phi_{i_2}\ev_3^*\phi_{i_3}}
{(\hb_1\!-\!\psi_1)(\hb_2\!-\!\psi_2)},
\end{equation}
once one computes
$$\wt\cZ_p(\hb_1,\hb_2,\al_{i_1},\al_{i_2},u)\equiv \sum_{d=0}^{\i}u^d
\int_{\ov\M_{0,3}(\P^{n-1},d)}\frac{\E(\V_0)\ev_1^*\phi_{i_1}\ev_2^*\phi_{i_2}\ev_3^*x^p}
{(\hb_1\!-\!\psi_1)(\hb_2\!-\!\psi_2)}.$$
Similarly, $\hb_1\hb_2\wt\cZ_p$ can be reconstructed from
the its ``mod $\hb_2^{-1}$'' part, if one can compute
$$\cZ_{pq}(\hb_1,\al_{i_1},u)\equiv \sum_{d=0}^{\i}u^d
\int_{\ov\M_{0,3}(\P^{n-1},d)}\frac{\E(\V_0)\ev_1^*\phi_{i_1}\ev_2^*x^p\ev_3^*x^q}
{\hb_1\!-\!\psi_1}.$$
Unfortunately, the ``mod $\hb_1^{-1}$ part'' of $\cZ_{pq}(\hb_1,x,u)$ is not simple.\\

\noindent
The above approach does, nevertheless, lead to a reconstruction theorem for
$\E(\V_0)$-twisted GW-invariants of~$\P^{n-1}$, or 
equivalently for GW-invariants of projective hypersurfaces.
The theorem arising here is different from~\cite{LeP} and~\cite{P} for example,
as the reduction is made to GW-invariants with low powers of $\psi$-classes
and without increasing the number of marked points.
Furthermore, in may be possible to get a handle on the ``components'' of the ``mod $\hb_1^{-1}$ part''  of $\cZ_{pq}$, i.e.
$$\int_{\ov\M_{0,3}(\P^{n-1},d)}\ev_1^*x^p\ev_2^*x^q\ev_3^*x^r,$$
and $\cZ_{pq}$ itself through the approach of
Subsections~\ref{bcov1-Zreg_subs} and~\ref{bcov1-localprppf_subs} in~\cite{bcov1}.
Along with~\e_ref{recur_dfn_e}, this approach leads to an explicit, but complicated,
recursion for $\cZ_{pq}$ or (the components of its ``mod $\hb_1^{-1}$ part'').

\section{Algebraic Observations} 
\label{alg_sec}

\subsection{On Rigidity of Certain Polynomial Conditions}
\label{genalb_subs}

\noindent
This subsection describes the extent of rigidity of power series with coefficients 
in rational functions that satisfy a certain recursion and a polynomiality condition.
It is the analogue of Section 30.3.\\

\noindent
Denote by $\bar\Z^+$ the set of nonnegative integers and by $[n]$, whenever $n\!\in\!\bar\Z^+$,
the set of positive number not exceeding~$n$:
$$\bar\Z^+=\big\{0,1,2,\ldots,\big\}, \qquad [n]=\big\{1,2,\ldots,n\big\}.$$
Whenever $f$ is a function of $w$ (and possibly of other variables)
which is holomorphic at $w\!=\!0$ (for a dense subspace of the other variables)
and $s\!\in\!\bar\Z^+$, let 
\begin{equation}\label{derivdfn_e}
\cD_w^sf=\frac{1}{s!}\bigg\{\frac{d}{dw}\bigg\}^s\!f(w)\Big|_{w=0}.
\end{equation}
This is a function of the other variables if there are any.\\

\noindent
Let $\ti\Q_{\al}$ be any field extension of $\Q_{\al}$, possibly $\Q_{\al}$ itself.
Given
$$Y\!\equiv\!Y(\hb,x,u),~Z\!\equiv\!Z(\hb,x,u)
\in \ti\Q_{\al}(\hb,x)\big[\big[u\big]\big],$$
we define
\begin{gather}
\Phi_{Y,Z}\!\equiv\!\Phi_{Y,Z}(\hb,u,z)\in\ti\Q_{\al}(\hb)\big[\big[u,z\big]\big] 
\qquad\hbox{by}\notag\\
\label{Phidfn_e}
\Phi_{Y,Z}(\hb,u,z)=
\sum_{i=1}^{i=n}\frac{e^{\al_iz}}{\prod_{k\neq i}(\al_i\!-\!\al_k)}
Y\big(\hb,\al_i,ue^{\hb z}\big)Z(-\hb,\al_i,u).
\end{gather}

\begin{lmm}
\label{Phistr_lmm1}
If $Y,Z\!\in\!\ti\Q_{\al}(\hb,x)[[u]]$, there exists a unique collection $$\big(E_{Y,Z;d}\!\equiv\!E_{Y,Z;d}(\hb,\Om)\big)_{d\in\bar\Z^+}
\subset\ti\Q_{\al}(\hb)[\Om]$$ 
such that the $\Om$-degree of $E_{Y,Z;d}$ is at most $(d\!+\!1)n\!-\!1$ 
for every $d\in\bar\Z^+$ and
\begin{equation}\label{phistr_lmm1e}
\Phi_{Y,Z}(\hb,u,z)= \sum_{d=0}^{\i}u^d\bigg( \frac{1}{2\pi\I}\oint e^{\Om z} 
\frac{E_{Y,Z;d}(\hb,\Om)}{\prod_{k=1}^{k=n}\prod_{r=0}^{r=d}(\Om\!-\!\al_k\!-\!r\hb)}d\Om
\bigg),
\end{equation}
where each path integral is taken over a simple closed loop in $\C$ enclosing all 
points  $\Om\!=\!\al_k\!+\!r\hb$ with $k\!=\!1,\ldots,n$ and $r\!=\!0,1,\ldots,d$.
The equality holds for a dense collection of complex parameters~$\hb$.
\end{lmm}

\noindent
{\it Proof:}
It can be assumed that
$$\al_k\!+\!r\hb \neq \al_{k'}\!+\!r'\hb \qquad\forall~ 
k,k'\in[n], ~r,r'\!\in\!\bar\Z^+, ~(r,k)\!\neq\!(r',k').$$
Note that for every $i\!=\!1,\ldots,n$ and $d'\!=\!0,1,\ldots,d$,
\begin{equation*}\begin{split}
&\prod_{r=0}^{r=d'-1}\!\!\!\!(\al_i\!+\!d'\hb\!-\!\al_i\!-\!r\hb)
\prod_{r=d'+1}^{r=d}\!\!\!\!(\al_i\!+\!d'\hb\!-\!\al_i\!-\!r\hb)~
\prod_{r=0}^{r=d}\prod_{k\neq i}(\al_i\!+\!d'\hb\!-\!\al_k\!-\!r\hb)\\
&\qquad 
=d'!\hb^{d'}(d\!-\!d')!(-\hb)^{d-d'} 
\bigg(\prod_{r=1}^{r=d'}\prod_{k\neq i}(\al_i\!-\!\al_k\!+\!r\hb)\bigg)
\bigg(\prod_{k\neq i}(\al_i\!-\!\al_k)\bigg)
\bigg(\prod_{r=1}^{r=d-d'}\prod_{k\neq i}(\al_i\!-\!\al_k\!-\!r\hb)\bigg)\\
&\qquad = \bigg(\prod_{k\neq i}(\al_i\!-\!\al_k)\bigg)
Q_{d'}(\hb,\al_i)\,Q_{d-d'}(-\hb,\al_i),
\end{split}\end{equation*}
where
\begin{equation}\label{phistr_lmm1e1a}
Q_d(\hb,x)\equiv \prod_{r=1}^{r=d}\prod_{k=1}^{k=n}(x\!-\!\al_k\!+\!r\hb)
\qquad\forall~d\!\in\!\bar\Z^+.
\end{equation}
By the Residue Theorem,
\begin{equation}\label{phistr_lmm1e1b}\begin{split}
&\frac{1}{2\pi\I}\oint e^{\Om z} 
\frac{E_{Y,Z;d}(\hb,\Om)}{\prod_{r=0}^{r=d}\prod_{k=1}^{k=n}(\Om\!-\!\al_k\!-\!r\hb)}d\Om\\
&\qquad=\sum_{d'=0}^{d'=d}\sum_{i=1}^{i=n}e^{(\al_i+d'\hb)z}
\frac{E_{Y,Z;d}(\hb,\al_i\!+\!d'\hb)}{\big(\prod_{k\neq  i}(\al_i\!-\!\al_k)\big)
Q_{d'}(\hb,\al_i)\,Q_{d-d'}(-\hb,\al_i)}\\
&\qquad=\sum_{d'=0}^{d'=d}\sum_{i=1}^{i=n}
\bigg(\frac{e^{\al_iz}}{\prod_{k\neq i}(\al_i\!-\!\al_k)}\bigg)
\bigg(\frac{(e^{\hb z})^{d'}}{Q_{d'}(\hb,\al_i)Q_{d-d'}(-\hb,\al_i)}\bigg)
E_{Y,Z;d}(\hb,\al_i\!+\!d'\hb).
\end{split}\end{equation}
On the other hand, since $Y,Z\!\in\!\ti\Q_{\al}(\hb,x)[[u]]$,
\begin{equation}\label{phistr_lmm1e3}
Y(\hb,x,u)=\sum_{d=0}^{\i}\frac{N_{Y;d}(\hb,\al_i)}{Q_d(\hb,\al_i)}u^d 
\quad\hbox{and}\quad
Z(\hb,x,u)=\sum_{d=0}^{\i}\frac{N_{Z;d}(\hb,\al_i)}{Q_d(\hb,\al_i)}u^d
\end{equation}
for unique $N_{Y;d},N_{Z;d}\!\in\!\ti\Q_{\al}(\hb,x)$.
By~\e_ref{Phidfn_e} and~\e_ref{phistr_lmm1e3},
\begin{equation}\label{phistr_lmm1e5}\begin{split}
\Phi_{Y,Z}(\hb,u,z)
&=\sum_{d=0}^{\i}\sum_{d'=0}^{d'=d}\sum_{i=1}^{i=n}
\frac{e^{\al_iz}}{\prod_{k\neq i}(\al_i\!-\!\al_k)}
\bigg(\frac{N_{Y;d'}(\hb,\al_i)}{Q_{d'}(\hb,\al_i)}\bigg)(ue^{\hb z})^{d'}
\bigg(\frac{N_{Z;d-d'}(-\hb,\al_i)}{Q_{d-d'}(-\hb,\al_i)}\bigg)u^{d-d'}\\
&=\sum_{d=0}^{\i}u^d\Bigg(\sum_{d'=0}^{d'=d}\sum_{i=1}^{i=n}
\frac{e^{\al_iz}}{\prod_{k\neq i}(\al_i\!-\!\al_k)}
\bigg(\frac{(e^{\hb z})^{d'}}{Q_{d'}(\hb,\al_i)Q_{d-d'}(-\hb,\al_i)}\bigg)\\
&\qquad\qquad\qquad\qquad\qquad\qquad\qquad\qquad
\times N_{Y;d'}(\hb,\al_i)N_{Z;d-d'}(-\hb,\al_i)\Bigg).
\end{split}\end{equation}
By \e_ref{phistr_lmm1e1b} and~\e_ref{phistr_lmm1e5}, \e_ref{phistr_lmm1e}
is satisfied if and only if
\begin{equation}\label{phistr_lmm1e7}
E_{Y,Z;d}(\hb,\al_i\!+\!d'\hb)
=N_{Y;d'}(\hb,\al_i)\cdot N_{Z;d-d'}(-\hb,\al_i)
\quad\forall\,i=1,\ldots,n,~d'=0,\ldots,d.
\end{equation}
For a dense collection of complex parameters $\hb$,
there exists a unique polynomial 
$$E_{Y,Z;d}(\hb,\Om)\in\ti\Q_{\al}(\hb)[\Om]$$ 
of $\Om$-degree at most $(d\!+\!1)n\!-\!1$ that satisfies~\e_ref{phistr_lmm1e7}.

\begin{lmm}
\label{Phistr_lmm2}
If $Y,Z\!\in\!\ti\Q_{\al}(\hb,x)[[u]]$ and
$(E_{Y,Z;d})_{d\in\bar\Z^+}\subset\ti\Q_{\al}(\hb)[\Om]$ are as in Lemma~\ref{Phistr_lmm1},
then
\begin{equation}\label{Phistr_lmm2e}\begin{split}
\Phi_{Y,Z}\in \ti\Q_{\al}[\hb]\big[\big[u,z\big]\big] 
\quad&\Llra\quad E_{Y,Z;d}\in\ti\Q_{\al}[\hb,\Om]~~\forall\, d\!\in\!\bar\Z^+\\
&\Llra\quad E_{Z,Y;d}\in\ti\Q_{\al}[\hb,\Om]~~\forall\, d\!\in\!\bar\Z^+ 
\quad\Llra\quad \Phi_{Z,Y}\in \ti\Q_{\al}[\hb]\big[\big[u,z\big]\big].
\end{split}\end{equation}\\
\end{lmm}

\noindent
{\it Proof:} 
The equivalence of the two middle statements in~\e_ref{Phistr_lmm2e}
follows from~\e_ref{phistr_lmm1e7}, which implies that
$$E_{Z,Y;d}(\hb,\Om)=E_{Y,Z;d}(-\hb,\Om\!-\!d\hb).$$
On the other hand, by the Residue Theorem on $S^2$,
\begin{equation}\label{phistr_lmm2e1}
\frac{1}{2\pi\I}\oint 
\frac{\Om^k d\Om}{\prod_{r=0}^{r=d}\prod_{k=1}^{k=n}(\Om\!-\!\al_k\!-\!r\hb)}
=\begin{cases}
0,&\hbox{if}~k<(d\!+\!1)n\!-\!1;\\
1,&\hbox{if}~k=(d\!+\!1)n\!-\!1;\\
R_{k-(d+1)n+1}^d(\hb),&\hbox{if}~k>(d\!+\!1)n\!-\!1,\\
\end{cases}
\end{equation}
where $R_s^d\!\in\!\Q_{\al}[\hb]$ is given by
$$R_s^d(\hb) = \cD_w^s
\bigg(\frac{1}{\prod_{r=0}^{r=d}\prod_{k=1}^{k=n}(1\!-\!(\al_k\!+\!r\hb)w)}\bigg)
\qquad\forall~ s\in\bar\Z^+.$$
The path integral in~\e_ref{phistr_lmm2e1} is again taken over a simple closed loop
enclosing all points $\Om\!=\!\al_k\!+\!r\hb$ with $r\!\le\!d$.
Write 
\begin{equation}\label{phistr_lmm2e3}
\Phi_{Y,Z}(\hb,u,z)=\sum_{d=0}^{\i}\sum_{q=0}^{\i}\frac{1}{q!}F_{d,q}(\hb)z^qu^d
\qquad\hbox{and}\qquad
E_{Y,Z;d}(\hb,\Om)=\sum_{s=0}^{(d+1)n-1}\!\!\!\!\!\!\!f_{d,s}(\hb)\Om^s.
\end{equation}
By \e_ref{phistr_lmm1e}, \e_ref{phistr_lmm2e1}, and~\e_ref{phistr_lmm2e3}, 
\begin{equation}\label{phistr_lmm2e5}\begin{split}
F_{d,q}(\hb) &=\sum_{s=0}^{(d+1)n-1}\!\!\!
\frac{1}{2\pi\I}\oint 
\frac{f_{d,s}(\hb)\,\Om^{q+s} d\Om}{\prod_{r=0}^{r=d}\prod_{k=1}^{k=n}(\Om\!-\!\al_k\!-\!r\hb)}\\
&=\sum_{s=\max(0,(d+1)n-1-q)}^{(d+1)n-1}\!\!\!\!\!\!\!\!\!\!\!\!\!\!\!
R_{q+s-(d+1)n+1}^d(\hb)\,f_{d,s}(\al;\hb).
\end{split}\end{equation}
Since $R_s^d\!\in\!\Q_{\al}[\hb]$, 
it follows that $F_{d,q}\!\in\!\ti\Q_{\al}[\hb]$ if $f_{d,s}\!\in\!\ti\Q_{\al}[\hb]$ for all~$s$.
Conversely, since $R_0^d(\hb)\!=\!1$,
$$F_{d,0},\ldots,F_{d,(d+1)n-1} \in \ti\Q_{\al}[\hb] \qquad\Lra\qquad
f_{d,(d+1)n-1},\ldots, f_{d,0}\in \ti\Q_{\al}[\hb].$$
These observations imply the two remaining statements of Lemma~\ref{Phistr_lmm2}.

\begin{dfn}
\label{recur_dfn}
For $i,j\!\in\![n]$ with $i\!\neq\!j$ and $d\!\in\!\Z^+$, let
$$C_i^j(d)=\frac{\prod_{r=0}^{nd-1}(n\al_i+r(\al_j\!-\!\al_i)/d\big)}
{d\,\underset{(r,k)\neq(d,j)}{\prod_{r=1}^{r=d}\prod_{k=1}^{k=n}}
(\al_i\!-\!\al_k+r(\al_j\!-\!\al_i)/d\big)}
\in \Q_{\al}.$$
We will call $Z\!\in\!\ti\Q_{\al}(\hb,x)[[u]]$ {\tt $C$-recursive} if 
\begin{equation}\label{recur_dfn_e}
Z(\hb,\al_i,u)=\sum_{d=0}^{\i}\bigg(\sum_{r=-N_d}^{r=N_d}\!\!\!Z_{i;d}^r\hb^{-r}\bigg)u^d
+\sum_{d=1}^{\i}\sum_{j\neq i}\frac{1}{\hb-\frac{\al_j-\al_i}{d}}
C_i^j(d)u^d Z\big((\al_j\!-\!\al_i)/d,\al_j,u\big)
\end{equation}
for every $i\!\in\![n]$ and for some $N_d\!\in\!\Z$ and
$Z_{i;d}^r\!\in\!\ti\Q_{\al}$.\footnote{The recursion (30.11) in~\cite{MirSym}
is a renormalization of the recursion~\e_ref{recur_dfn_e} 
with a slightly different coefficient~$C_i^j(d)$.} 
\end{dfn}

\begin{prp}
\label{uniqueness_prp}
Suppose 
$$Y,Z\in\ti\Q_{\al}(\hb,x)[[u]]\subset 
\ti\Q_{\al}(x)\big[\big[\hb^{-1},u\big]\big]
+\ti\Q_{\al}(x)[\hb]\big[\big[u\big]\big]$$
are such that $Z$ is $C$-recursive, $\Phi_{Y,Z}\!\in\!\ti\Q_{\al}[\hb][[u,z]]$, 
and for every $i\!\in\![n]$
$$Y(\hb,\al_i,f_i) \equiv f_i ~~(\mod u)$$
for some $f_i\!\in\!\ti\Q_{\al}^*$. Then,
$$Z(\hb,\al_i,u)\equiv 0 ~~(\mod \hb^{-1}) ~~~\forall\,i\in[n]
\qquad\Lra\qquad Z(\hb,\al_i,u)= 0 ~~~\forall\,i\in[n].$$\\
\end{prp}

\noindent
{\it Remark 1:} Suppose 
$$Z(\hb,\al_i,u)=\sum_{d=0}^{\i}
\bigg(\sum_{r=-N_d}^{\i}\!\!\!\ti{Z}_{i;d}^r\hb^{-r}\bigg)u^d$$
for some $\ti{Z}_{i;d}^r\!\in\!\ti\Q_{\al}$.
In the statement of Proposition~\ref{uniqueness_prp} and throughout the rest of the paper, 
$$Z(\hb,\al_i,u) \equiv \sum_{d=0}^{\i}
\bigg(\sum_{r=-N_d}^{a-1}\!\!\!\ti{Z}_{i;d}^r\hb^{-r}\bigg)u^d
~~(\mod \hb^{-a}),$$ 
i.e.~we drop $\hb^{-a}$ and higher powers of $\hb^{-1}$, 
instead of higher powers of~$\hb$.\\

\noindent
{\it Remark 2:} In contrast to the situation in \cite[Chapter~30]{MirSym},
the assumptions of Proposition~\ref{uniqueness_prp}, 
i.e.~recursivity and $\hb$-polynomiality with respect to~$Y$, are both linear conditions
on~$Z$. 
Consequently, the $\hb^{-1}$-term of $Z(\hb,\al_i,u)$ is no longer necessary to 
determine~$Z$.\\

\noindent
{\it Proof:} Suppose $d\!\ge\!0$ and 
we have shown that
\begin{equation}\label{polyn_prp_e1}
Z(\hb,\al_i,u)\equiv 0 ~~(\mod u^d) \qquad\forall~i=1,\ldots,n.
\end{equation}
With notation as in~\e_ref{recur_dfn_e} and by the last assumption on $Z(\hb,\al_i,u)$, 
it follows that
\begin{equation}\label{polyn_prp_e2}
Z(\hb,\al_i,u)\equiv u^d\sum_{r=1}^{r=N_d}\!\!\!Z_{i;d}^r\hb^{-r}
~~(\mod u^{d+1}) \qquad\forall~i=1,\ldots,n.
\end{equation}
If $N_{Y;d'},N_{Z;d'}\!\in\!\ti\Q_{\al}(\hb,x)$ are as in the proof of 
Lemma~\ref{Phistr_lmm1},
\begin{equation}\label{polyn_prp_e3}
N_{Y;0}(\hb,\al_i)=f_i \qquad\hbox{and}\qquad
N_{Z;d'}(\hb,\al_i)=\begin{cases}
0,&\hbox{if}~d'\!<\!d;\\
Q_d(\hb,\al_i)\sum_{r=1}^{r=N_d}Z_{i;d}^r\hb^{-r},&\hbox{if}~d'\!=\!d,
\end{cases}\end{equation}
by \e_ref{phistr_lmm1e3}, \e_ref{polyn_prp_e1}, and \e_ref{polyn_prp_e2}.
Since
$$E_{Z,Y;d}(\hb,\al_i\!+\!d'\hb)=0 ~~\forall\, 
d'=0,1,\ldots,d\!-\!1,~i=1,\ldots,n$$
by~\e_ref{phistr_lmm1e7} and~\e_ref{polyn_prp_e3} and 
$E_{Z,Y;d}\!\in\!\ti\Q_{\al}[\hb,\Om]$ by Lemma~\ref{Phistr_lmm2},
$$E_{Z,Y;d}(\hb,\Om)=\bigg(\prod_{d'=0}^{d-1}\prod_{j=1}^{j=n}
\big(\Om\!-\!\al_j\!-\!d'\hb\big)\bigg)\cdot R_d(\hb,\Om)$$
for some $R_d\!\in\!\ti\Q_{\al}[\hb,\Om]$.
Thus,
\begin{equation}\label{polyn_prp_e5}
E_{Z,Y;d}(\hb,\al_i\!+\!d\hb)
=\bigg(\prod_{d'=0}^{d-1}\prod_{k=1}^{k=n}
\big((\al_i\!+\!d\hb)\!-\!\al_k\!-\!d'\hb\big)\bigg)\cdot R_d(\hb,\al_i\!+\!d\hb)
=\hb^d\ti{R}_d(\hb)
\end{equation}
for some $\ti{R}_d\!\in\!\ti\Q_{\al}[\hb]$.
On the other hand, by~\e_ref{phistr_lmm1e7} and~\e_ref{polyn_prp_e3} 
\begin{equation}\label{polyn_prp_e6}\begin{split}
E_{Z,Y;d}(\hb,\al_i\!+\!d\hb) =N_{Z;d}(\hb,\al_i)\cdot f_i
&=f_i\cdot \bigg(d!\hb^d\prod_{r=1}^{r=d}\prod_{k\neq i}(\al_i\!-\!\al_k\!+\!r\hb)
\bigg) \sum_{r=1}^{r=N_d}\!\!Z_{i;d}^r\hb^{-r}\\
&=f_i\cdot
\bigg(d!\prod_{r=1}^{r=d}\prod_{k\neq i}(\al_i\!-\!\al_k\!+\!r\hb)\bigg)
\sum_{r=1}^{r=N_d}\!\!Z_{i;d}^r\hb^{d-r}.
\end{split}\end{equation}
By~\e_ref{polyn_prp_e5} and~\e_ref{polyn_prp_e6},
$$Z_{i;d}^r=0 ~~~\forall~ r=1,\ldots,N_d,~i=1,\ldots,n.$$
Along with~\e_ref{polyn_prp_e2}, this implies that \e_ref{polyn_prp_e1}
holds with $d$ replaced by~$d\!+\!1$.

\subsection{Admissible Transforms}
\label{tranform_subs}

\noindent
This subsection is the analogue of the beginning of Section 30.4 in~\cite{MirSym}.
We describe transforms of $Y,Z\!\in\!\ti\Q_{\al}(\hb,x)$ that preserve the polynomiality 
property of Lemma~\ref{Phistr_lmm2} and the recursivity property of Definition~\ref{recur_dfn}.
The statement of Lemma~\ref{Phistr_lmm4} is followed by complete proofs.
The first of the five transforms below has no analogue in~\cite{MirSym}.

\begin{lmm}
\label{Phistr_lmm4}
Suppose $Y,Z\in\ti\Q_{\al}(\hb,x)[[u]]$ are such that $Z$ is $C$-recursive and $\Phi_{Y,Z}\!\in\!\ti\Q_{\al}[\hb][[u,z]]$.
Then, 
\begin{enumerate}[label=(\roman*)]
\item\label{deriv_ch} if $u\!=\!e^t$, 
$\bar{Z}\!\equiv\!\big\{x\!+\!\hb\frac{d}{dt}\big\}Z$
is $C$-recursive and $\Phi_{Y,\bar{Z}}\!\in\!\ti\Q_{\al}[\hb][[u,z]]$;

\item\label{mult_ch} if $f\!\in\!\ti\Q_{\al}[u]$, then 
$fZ$ is $C$-recursive and $\Phi_{Y,fZ},\Phi_{fY,Z}\!\in\!\ti\Q_{\al}[\hb][[u,z]]$;

\item\label{hbmult_ch} if $f\!\in\!\ti\Q_{\al}[\hb]$, then 
$\bar{Z}\!\equiv\!fZ$ is $C$-recursive and $\Phi_{Y,\bar{Z}}\!\in\!\ti\Q_{\al}[\hb][[u,z]]$;

\item\label{multoverh_ch} if $f\!\in\!\ti\Q_{\al}[u]$, $f(0)\!=\!0$,
and $\bar{Y}\!=\!e^{f/\hb}Y$, then
$\bar{Z}\!\equiv\!e^{f/\hb}Z$ is $C$-recursive and 
$\Phi_{\bar{Y},\bar{Z}}\!\in\!\ti\Q_{\al}[\hb][[u,z]]$;

\item\label{transf_ch} if $g\!\in\!\Q[[u]]$, $g(0)\!=\!0$, 
$$\bar{Y}(\hb,x,u)=e^{x g(u)/\hb}Y\big(\hb,x,u e^{g(u)}\big), \qquad\hbox{and}\qquad
\bar{Z}(\hb,x,u)=e^{x g(u)/\hb}Z\big(\hb,x,u e^{g(u)}\big),$$
then $\bar{Z}$ is $C$-recursive and  $\Phi_{\bar{Y},\bar{Z}}\!\in\!\ti\Q_{\al}[\hb][[u,z]]$.\\
\end{enumerate}
\end{lmm}

\noindent
{\it Remark:} In fact, \ref{mult_ch} and \ref{hbmult_ch} are special cases of
the admissible transform defined by $f\!\in\!\Q_{\al}[\hb][[u]]$.\\

\noindent
\ref{deriv_ch} The operator $\big\{\al_i\!+\!\hb\frac{d}{dt}\big\}$ preserves
the structure of the first term on the right-hand side of~\e_ref{recur_dfn_e}.
The $(d,j)$-summand in the last term becomes
\begin{alignat}{1}
\Big\{\al_i\!+\!\hb\frac{d}{dt}\Big\} &
\bigg( \frac{C_i^j(d)u^d}{\hb-\frac{\al_j-\al_i}{d}} 
Z\big((\al_j\!-\!\al_i)/d,\al_j,u\big) \bigg)
=\frac{C_i^j(d)u^d}{\hb-\frac{\al_j-\al_i}{d}} 
\Big\{\al_i\!+\!d\hb\!+\!\hb\frac{d}{dt}\Big\}Z\big((\al_j\!-\!\al_i)/d,\al_j,u\big)\notag\\
\label{deriv_ch_e1}
&\qquad\qquad
=\frac{C_i^j(d)u^d}{\hb-\frac{\al_j-\al_i}{d}}\bar{Z}\big((\al_j\!-\!\al_i)/d,\al_j,u\big)
+dC_i^j(d)u^dZ\big((\al_j\!-\!\al_i)/d,\al_j,u\big).
\end{alignat}
Since the last term in~\e_ref{deriv_ch_e1} does not depend on $\hb$ and 
$Z$ is $C$-recursive, it follows that $\bar{Z}$ is also $C$-recursive.
Since 
$$\Phi_{\bar{Z},Y}(\hb,u,z)=\frac{d}{dz}\Phi_{Z,Y}(\hb,u,z)$$
and $\Phi_{Y,Z}\!\in\!\ti\Q_{\al}[\hb][[u,z]]$, 
$\Phi_{Y,\bar{Z}}\!\in\!\ti\Q_{\al}[\hb][[u,z]]$ by the middle equivalence
in Lemma~\ref{Phistr_lmm2}.\\

\noindent
\ref{mult_ch} Since $Z$ is $C$-recursive and the multiplication by $f$
preserves the structure of each of the terms in~\e_ref{recur_dfn_e},
$fZ$ is also $C$-recursive.
Since $\Phi_{Y,fZ}\!=\!f\Phi_{Y,Z}$ and $\Phi_{Y,Z}\!\in\!\ti\Q_{\al}[\hb][[u,z]]$, 
$\Phi_{Y,fZ}\!\in\!\ti\Q_{\al}[\hb][[u,z]]$.
Similarly, since $\Phi_{Z,fY}\!=\!f\Phi_{Z,Y}$, 
$\Phi_{fY,Z}\!\in\!\ti\Q_{\al}[\hb][[u,z]]$ by Lemma~\ref{Phistr_lmm2}.\\

\noindent
\ref{hbmult_ch} It is sufficient to assume that $f(\hb)\!=\!\hb$.
The multiplication by $\hb$ preserves
the structure of the first term on the right-hand side of~\e_ref{recur_dfn_e}.
The $(d,j)$-summand in the last term becomes
\begin{equation}\label{hbmult_ch_e1}\begin{split}
\hb\frac{C_i^j(d)u^d}{\hb-\frac{\al_j-\al_i}{d}} 
Z\big((\al_j\!-\!\al_i)/d,\al_j,u\big) 
&=\frac{C_i^j(d)u^d}{\hb-\frac{\al_j-\al_i}{d}} 
\bar{Z}\big((\al_j\!-\!\al_i)/d,\al_j,u\big)\\
&\qquad +C_i^j(d)u^dZ\big((\al_j\!-\!\al_i)/d,\al_j,u\big).
\end{split}\end{equation}
Since the last term in~\e_ref{hbmult_ch_e1} does not depend on $\hb$ and 
$Z$ is $C$-recursive, it follows that $\bar{Z}$ is also $C$-recursive.
Since $\Phi_{Y,\bar{Z}}\!=\!-\hb\Phi_{Y,Z}$ and $\Phi_{Y,Z}\!\in\!\ti\Q_{\al}[\hb][[u,z]]$, 
$\Phi_{Y,\bar{Z}}\!\in\!\ti\Q_{\al}[\hb][[u,z]]$.\\

\noindent
\ref{multoverh_ch}
Since $f(0)\!=\!0$, i.e.~$f$ contains no degree-0 term in $u$,
the multiplication by $e^{f(u)/\hb}$ preserves
the structure of the first term on the right-hand side of~\e_ref{recur_dfn_e}.
The $(d,j)$-summand in the last term becomes
\begin{equation*}\begin{split}
&e^{f(u)/\hb}\frac{C_i^j(d)u^d}{\hb-\frac{\al_j-\al_i}{d}} 
Z\big((\al_j\!-\!\al_i)/d,\al_j,u\big)
=\frac{C_i^j(d)u^d}{\hb-\frac{\al_j-\al_i}{d}} 
\bar{Z}\big((\al_j\!-\!\al_i)/d,\al_j,u\big)\\
&\qquad\qquad\qquad\qquad\qquad+\Big(e^{f(u)/\hb}-e^{f(u)/((\al_j-\al_i)/d)}\Big)
\frac{C_i^j(d)u^d}{\hb-\frac{\al_j-\al_i}{d}}
Z\big((\al_j\!-\!\al_i)/d,\al_j,u\big).
\end{split}\end{equation*}
Since $Z$ is $C$-recursive and 
$$\frac{e^{f(u)/\hb}-e^{f(u)/((\al_j-\al_i)/d)}}{\hb-\frac{\al_j-\al_i}{d}}
 \in \ti\Q_{\al}[\hb,\hb^{-1}]\big[\big[u\big]\big],$$
it follows that $\bar{Z}$ is $C$-recursive as well.
On the other hand,
\begin{equation}\label{multoverh_ch_e2}\begin{split}
\Phi_{\bar{Y},\bar{Z}}(\hb,u,z)=e^{(f(ue^{\hb z})-f(u))/\hb}
\Phi_{Y,Z}(\hb,u,z).
\end{split}\end{equation}
Since $\Phi_{Y,Z}\!\in\!\ti\Q_{\al}[\hb][[u,z]]$ and
$$\big(f(ue^{\hb z})-f(u)\big)/\hb \in \ti\Q_{\al}[\hb]\big[\big[u,z\big]\big],$$
\e_ref{multoverh_ch_e2} implies that 
$\Phi_{\bar{Y},\bar{Z}}\!\in\!\ti\Q_{\al}[\hb][[u,z]]$ as well.\\

\noindent
\ref{transf_ch}
Since $g(0)\!=\!0$, the operation of replacing $u$ with $ue^{g(u)}$ followed
by multiplication by $e^{\al_ig(u)/\hb}$ preserves the structure of the first term 
on the right-hand side of~\e_ref{recur_dfn_e}.
The $(d,j)$-summand in the last term becomes
\begin{equation*}\begin{split}
&e^{\al_i g(u)/\hb}\frac{C_i^j(d)u^de^{dg(u)}}{\hb-\frac{\al_j-\al_i}{d}} 
Z\big((\al_j\!-\!\al_i)/d,\al_j,u\big)
=\frac{C_i^j(d)u^d}{\hb-\frac{\al_j-\al_i}{d}} 
\bar{Z}\big((\al_j\!-\!\al_i)/d,\al_j,u\big)\\
&\qquad\qquad\qquad\qquad\qquad+
\Big(e^{(\al_i/\hb+d)g(u)}-e^{(\al_j/((\al_j-\al_i)/d))g(u)}\Big)
\frac{C_i^j(d)u^d}{\hb-\frac{\al_j-\al_i}{d}}
Z\big((\al_j\!-\!\al_i)/d,\al_j,u\big).
\end{split}\end{equation*}
Since $Z$ is $C$-recursive and 
$$\frac{e^{(\al_i/\hb+d)g(u)}-e^{(\al_j/((\al_j-\al_i)/d))g(u)}}
{\hb-\frac{\al_j-\al_i}{d}}
=\frac{d}{(\al_i\!+\!d\hb)-\al_j}e^{dz/(z-\al_i)}\Big|_{z=\al_j}^{z=\al_i+d\hb}
 \in \Q_{\al}[\hb,\hb^{-1}]\big[\big[u\big]\big],$$
it follows that $\bar{Z}$ is $C$-recursive as well.
On the other hand,
\begin{equation}\label{transf_ch_e2}\begin{split}
\Phi_{\bar{Y},\bar{Z}}(\hb,u,z)=\Phi_{Y,Z}\big(\hb,ue^{g(u)},\ti{z}\big),
\qquad\hbox{where}\quad
\ti{z}=z+\frac{g(ue^{\hb})\!-\!g(u)}{\hb}\in\Q[\hb,z]\big[\big[u\big]\big].
\end{split}\end{equation}
Since $\Phi_{Y,Z}\!\in\!\ti\Q_{\al}[\hb][[u,z]]$,
\e_ref{transf_ch_e2} implies that 
$\Phi_{\bar{Y},\bar{Z}}\!\in\!\ti\Q_{\al}[\hb][[u,z]]$ as well.

\subsection{Some Properties of Hypergeometric Series}
\label{algcomp_subs}

\noindent
In this subsection we verify the three claim concerning the power series $\cY_p(\hb,x,u)$
made in Subsection~\ref{outline_subs}:
\begin{enumerate}[label=(\alph*)]
\item\label{Ymprec_item} $\cY_{-1}(\hb,x,u)$ satisfies the $C$-recursivity condition
of Definition~\ref{recur_dfn};
\item\label{YmpPhi_item} $\Phi_{\cY,\cY_{-1}}\!\in\!\Q_{\al}[\hb][[u,z]]$;
\item\label{Ypcontr_item} there exist $\ti{C}_{p,q}^{(r)}$ as in Theorem~\ref{main_thm}
such that \e_ref{Ypcond_e} is satisfied.\\
\end{enumerate}

\noindent
{\it Proof of~\ref{Ymprec_item}:}
Note that 
$$\frac{C_i^j(d)u^d}{\hb-\frac{\al_j-\al_i}{d}} 
\, \cY_{-1}\big((\al_j\!-\!\al_i)/d,\al_j,u\big)
=\Res_{z=\frac{\al_j-\al_i}{d}}\bigg\{\frac{1}{\hb\!-\!z}\cY_{-1}(z,\al_i,u)\bigg\}.$$
Thus, by the Residue Theorem on~$S^2$,
\begin{equation}\label{recurpf_e1}\begin{split}
\sum_{d=1}^{\i}\sum_{j\neq i}\frac{C_i^j(d)u^d}{\hb-\frac{\al_j-\al_i}{d}}
&\cY_{-1}\big((\al_j\!-\!\al_i)/d,\al_j,u\big)
=-\Res_{z=\hb,0,\i}\bigg\{\frac{1}{\hb\!-\!z}\cY_{-1}(z,\al_i,u)\bigg\}\\
&\qquad\qquad\qquad\qquad
=\cY_{-1}(\hb,\al_i,u)-\Res_{z=0,\i}\bigg\{\frac{1}{\hb\!-\!z}\cY_{-1}(z,\al_i,u)\bigg\}.
\end{split}\end{equation}
On the other hand, if $\cY_{-1;d}$ is the degree-$d$ term of $\cY_{-1}$,
\begin{alignat*}{1}
\Res_{z=0}\bigg\{\frac{1}{\hb\!-\!z}\cY_{-1;d}(z,\al_i,u)\bigg\}
&=\cD^{d-1}_z\Bigg\{\frac{1}{\hb\!-\!z}\frac{\prod_{r=0}^{r=nd-1}(n\al_i\!+\!rz)}
{d!\prod_{r=1}^{r=d}\prod_{k\neq i}(\al_i\!-\!\al_k\!+\!rz)}\Bigg\}\in\Q_{\al}[\hb^{-1}],\\
\Res_{z=\i}\bigg\{\frac{1}{\hb\!-\!z}\cY_{-1}(z,\al_i,u)\bigg\}&=1,
\end{alignat*}
where $\cD^d_z$ is as in~\e_ref{derivdfn_e}.
Thus, \e_ref{recurpf_e1} implies that $\cY_{-1}$ satisfies the recursion~\e_ref{recur_dfn_e}.\\

\noindent
{\it Proof of~\ref{YmpPhi_item}:}
Let 
$$R(\hb,x,u)=\sum_{d=0}^{\i}u^d \frac{\prod_{r=1}^{r=nd}(nx\!+\!r\hb)}
{\prod_{r=1}^{r=d}\prod_{k=1}^{k=n}(x\!-\!\al_k\!+\!r\hb)}.$$
By \ref{mult_ch} of Lemma~\ref{Phistr_lmm4}, it is sufficient to show that $\Phi_{R,\cY_{-1}}\!\in\!\Q_{\al}[\hb][[u,z]]$.
Note that 
$$\frac{e^{\al_iz}}{\prod\limits_{k\neq i}(\al_i\!-\!\al_k)}
R\big(\hb,\al_i,ue^{\hb z}\big)\cY_{-1}(-\hb,\al_i,u)
=\Res_{x=\al_i}\Bigg\{\frac{e^{xz}}{\prod\limits_{k=1}\limits^{k=n}(x\!-\!\al_k)}
R\big(\hb,x,ue^{\hb z}\big)\cY_{-1}(-\hb,x,u)\Bigg\}.$$
Thus, by the Residue Theorem on~$S^2$,
\begin{equation*}\begin{split}
\Phi_{R,\cY_{-1}}(\hb,u,z) &=-
\Res_{x=\i}\Bigg\{\frac{e^{xz}}{\prod\limits_{k=1}\limits^{k=n}(x\!-\!\al_k)}
R\big(\hb,x,ue^{\hb z}\big)\cY_{-1}(-\hb,x,u)\Bigg\}\\
&=\sum_{p=0}^{\i}\frac{z^{n-1+p}}{(n\!-\!1\!+\!p)!}
\cD_w^p\Bigg\{\frac{1}{\prod\limits_{k\neq i}(1\!-\!\al_kw)}
\bigg(\sum_{d=0}^{\i}u^d e^{d\hb z}\frac{\prod_{r=1}^{r=nd}(n\!+\!r\hb w)}
{\prod_{r=1}^{r=d}\prod_{k=1}^{k=n}(1\!-\!(\al_k\!-\!r\hb)w)}\bigg)\\
&\qquad\qquad\qquad\qquad\qquad\qquad\qquad\qquad
\times\bigg(\sum_{d=0}^{\i}u^d \frac{\prod_{r=0}^{r=nd-1}(n\!-\!r\hb w)}
{\prod_{r=1}^{r=d}\prod_{k=1}^{k=n}(1\!-\!(\al_k\!+\!r\hb)w)}\bigg)\Bigg\}.
\end{split}\end{equation*}
The $p$th summand above is polynomial in $\hb$. Thus, 
$\Phi_{R,\cY_{-1}}\!\in\!\Q_{\al}[\hb][[u,z]]$.\\

\noindent
{\it Proof of \ref{Ypcontr_item}:}
Suppose $p\!\in\!\Z^+$, and we have constructed power series 
$$\cY_0,\cY_1,\ldots,\cY_{p-1}\in\Q_{\al}(\hb,x)\big[\big[e^t\big]\big],$$
satisfying the second equation in~\e_ref{main_e1} and \e_ref{Ypcond_e}, so that  
\begin{equation}\label{Yexpand_e5}
e^{xt/\hb}\,\cY_{p-1}(\hbar,x,e^t)=x^p+
\sum_{q=1}^{\i}\bigg(\sum_{r=0}^{p-1+q}C_{p-1,q}^{(r)}(t)\,x^{p+q-r}\bigg)\hb^{-q},
\end{equation}
where $C_{p-1,q}^{(r)}(t)$  is a degree-$r$ symmetric polynomial in $\al_1,\ldots,\al_n$
with coefficients in $\Q[t][[e^t]]$  such that
\begin{equation}\label{Yexpand_e6}
C_{p-1,q}^{(0)}(t)=I_{p-1,p-1+q}(t)\big/I_{p-1,p-1}(t).
\end{equation}
This assumption is satisfied for $p\!=\!1$ by~\e_ref{Yexpand_e} and~\e_ref{Yexpand_e2}. 
Since $\cY_{p-1}\in\Q_{\al}(\hb,x)[[e^t]]$,
\begin{gather}
C_{p-1,1}^{(0)}(t)\in t+\Q_{\al}[[e^t]],\quad
C_{p-1,1}^{(r)}(t)\in \Q_{\al}[[e^t]] ~~\forall\,r\!\ge\!1
\qquad\Lra\notag\\
\label{psstr_e}
\frac{d}{dt}C_{p-1,1}^{(r)}(t)\in \Q_{\al}[[e^t]] ~~\forall\,r\!\ge\!0.
\end{gather}
Thus,
\begin{equation}\label{Yexpand_e8}\begin{split}
\cY_p(\hb,x,e^t)&\equiv e^{-xt/\hb}\frac{\hb}{I_{p,p}(t)}\frac{d}{dt} \big(e^{xt/\hb}\cY_{p-1}(\hb,x,e^t)\big)\\
&\qquad\qquad
-\frac{1}{I_{p,p}(t)}\sum_{r=1}^{r=p} \Big(\frac{d}{dt}C_{p-1,1}^{(r)}(t)\Big)\cY_{p-r}(\hb,x,e^t)
\in\Q_{\al}(\hb,x)\big[[e^t]\big].
\end{split}\end{equation}
By~\e_ref{Tden_e}, \e_ref{Yexpand_e5}, and \e_ref{Yexpand_e6} are
satisfied with $p$ replaced by $p\!+\!1$ and
$$C_{p,q}^{(r)}(t) = \frac{1}{I_{p,p}(t)}\frac{d}{dt}C_{p-1,q+1}^{(r)}(t)
-\frac{1}{I_{p,p}(t)}\sum_{s=1}^{\min(p,r)}\Big(\frac{d}{dt} C_{p-1,1}^{(s)}(t)\Big)
C_{p-s,q}^{(r-s)}(t).$$
In particular, $\cY_p$ satisfies \e_ref{Ypcond_e}.
By~\e_ref{Yexpand_e8}, the coefficients $\ti{C}^{(r)}_{p,q}$ are inductively defined~by
\begin{gather}\label{coeffind_e}
\ti{C}^{(r)}_{p,q}(e^t)=I_{p,p}(t)^{-1}\Bigg(
\frac{d}{dt}\ti{C}_{p-1,q}^{(r)}+
I_{q,q}(t)\ti{C}_{p-1,q-1}^{(r)}(e^t)
-\sum_{s=1}^{r-1}\bigg(\frac{d}{dt}C_{p-1,1}^{(s)}\bigg)\ti{C}_{p-s,q}^{(q-s)}\Bigg),\\
\hbox{where}\qquad \ti{C}_{p-1,p-r}^{(r)}\equiv -C_{p-1,1}^{(r)}, 
\qquad \ti{C}_{p-1,-1}^{(r)}\equiv0, 
\qquad \ti{C}_{p',q}^{(r)}\equiv0~~\forall\,r\!\le\!0.\notag
\end{gather}
Thus, by~\e_ref{psstr_e} and inductive assumptions, $\ti{C}^{(r)}_{p,q}\!\in\!\Q_{\al}[[e^t]]$
is a degree-$r$ symmetric polynomial in $\al_1,\ldots,\al_n$ with coefficients in 
$\Q_{\al}[[e^t]]$, as required.

\section{Localization Computations}
\label{localization_sec}

\subsection{Equivariant Cohomology}
\label{equivcoh_subs}

\noindent
In Subsection~\ref{local_subs1}, we apply the classical localization theorem~\e_ref{ABothm_e}
with the standard action of the $n$-torus $\T$ on $\P^{n-1}$ and on
$\ov\M_{0,m}(\P^{n-1},d)$ to verify Lemma~\ref{recgen_lmm}.
In Subsection~\ref{local_subs2}, we apply~\e_ref{ABothm_e} with an action of
the $(n\!+\!1)$-torus
$$\wt\T\equiv \T\times\T^1$$
on $\P^{n-1}\!\times\!\P^1$ and on a subspace of 
$\ov\M_{0,m}(\P^1\!\times\!\P^{n-1},(1,d))$ to verify Lemma~\ref{polgen_lmm}.
The aim of this subsection is to review the basics of equivariant cohomology and
to set notation.
Throughout this subsection, $G$ denotes an $m$-torus, 
either $(\C^*)^m$ or~$(S^1)^m$.\\

\noindent
The $m$-torus $G$ acts freely on $EG\!=\!(\C^{\i})^m\!-\!0$ (or $(S^{\i})^m$):
$$\big(e^{\I\th_1},\ldots,e^{\I\th_m}\big)\cdot (z_1,\ldots,z_m)
=\big(e^{\I\th_1}z_1,\ldots,e^{\I\th_m}z_m\big).$$
Thus, the classifying space for $G$ and its group cohomology are given by
$$\qquad BG\equiv EG/G=(\P^{\i})^m  \qquad\hbox{and}\qquad
H_G^*\equiv H^*(BG;\Q)=\Q[\al_1,\ldots,\al_m],$$
where $\al_i\!=\!\pi_i^*c_1(\ga^*)$ if
$$\pi_i\!: (\P^{\i})^m\lra\P^{\i} \qquad\hbox{and}\qquad
\ga\lra\P^{\i}$$
are the projection onto the $i$th component and the tautological line bundle, respectively.
Denote by $\H_G^*$ the field of fractions of $H_G^*$:
$$\H_G^*=\Q_{\al}\equiv \Q(\al_1,\ldots,\al_m).$$\\

\noindent
A representation $\rho$ of $G$, i.e.~a linear action of $G$ on $\C^k$,
induces a vector bundle over~$BG$:
$$V_{\rho}\equiv EG\times_G\C^k.$$
If $\rho$ is one-dimensional, we will call
$$c_1(V_{\rho}^*)=-c_1(V_{\rho})\in H_G^*\subset \H_G^*$$ 
the {\tt weight} of $\rho$.
For example, $\al_i$ is the weight of representation
\begin{equation}\label{ifactorrep_e}
\pi_i\!: G\lra\C^*, \qquad 
\big(e^{\I\th_1},\ldots,e^{\I\th_m}\big)\cdot z =  e^{\I\th_i}z.
\end{equation}
More generally, if a representation $\rho$ of $G$ on $\C^k$ splits 
into one-dimensional representations with weights $\be_1,\ldots,\be_k$,
we will call $\be_1,\ldots,\be_k$ the {\tt weights} of~$\rho$.
In such a case,
\begin{equation}\label{weightspord_e}
e(V_{\rho}^*)=\be_1\cdot\ldots\cdot\be_k.
\end{equation}
We will call the representation $\rho$ of $\T$ on $\C^n$ with weights $\al_1,\ldots,\al_n$
the {\tt standard representation} of~$\T$.\\

\noindent
If $G$ acts on a topological space $M$, let
$$H_G^*(M)\equiv H^*(BG;\Q), \qquad\hbox{where}\qquad BM=EG\!\times_G\!M,$$
be the {\tt equivariant cohomology} of $M$.
The projection map $BM\!\lra\!BG$ induces an action of $H_G^*$ on $H_G^*(M)$.
Let
$$\H_G^*(M)=H_G^*(M)\otimes_{H_G^*}\H_G^*.$$
If the $G$-action on $M$ lifts to an action on a (complex) vector bundle $V\!\lra\!M$,
then 
$$BV\equiv EG\!\times_G\!V$$
is a vector bundle over $BM$.
Let
$$\E(V)\equiv e(BV)\in H_G^*(M)\subset \H_G^*(M)$$
denote the {\tt equivariant euler class of} $V$.\\

\noindent
The standard action of $\T$ on $\P^{n-1}$ is 
the action induced by the standard action $\rho$ of $\T$ on~$\C^n$:
$$\big(e^{\I\th_1},\ldots,e^{\I\th_n}\big)\cdot [z_1,\ldots,z_n] 
=\big[e^{\I\th_1}z_1,\ldots,e^{\I\th_n}z_n\big].$$
Since $B\P^{n-1}=\P V_{\rho}$, 
$$H_{\T}^*(\P^{n-1})\equiv H^*\big(\P V_{\rho};\Q\big)
= \Q[x,\al_1,\ldots,\al_n]
\big/\big(x^n\!+\!c_1(V_{\rho})x^{n-1}\!+\!\ldots\!+\!c_n(V_{\rho})\big),$$
where $x\!=\!c_1(\ti\ga^*)$ and $\ti\ga\!\lra\!\P V_{\rho}$ 
is the tautological line bundle.
Since 
$$c(V_{\rho})=(1-\al_1)\ldots(1-\al_n),$$
it follows that 
\begin{equation}\label{pncoh_e}\begin{split}
H_{\T}^*(\P^{n-1}) &= \Q[x,\al_1,\ldots,\al_n]\big/(x\!-\!\al_1)\ldots(x\!-\!\al_n)
\qquad\hbox{and}\\
\H_{\T}^*(\P^{n-1}) &= \Q_{\al}[x]\big/(x\!-\!\al_1)\ldots(x\!-\!\al_n).
\end{split}\end{equation}\\

\noindent
The standard action of $\T$ on $\P^{n-1}$ has $n$-fixed points:
$$P_1=[1,0,\ldots,0], \qquad P_2=[0,1,0,\ldots,0], 
\quad\ldots\quad P_n=[0,\ldots,0,1].$$
For each $i\!=\!1,2,\ldots,n$, let
\begin{equation}\label{phidfn_e}
\phi_i= \prod_{k\neq i}(x\!-\!\al_k) \in H_{\T}^*(\P^{n-1}).
\end{equation}
By equation~\e_ref{phiprop_e} below, $\phi_i$ is the equivariant Poincare dual of $P_i$.
We also note that $\ti\ga|_{BP_i}\!=\!V_{\pi_i}$, where $\pi_i$ is as in~\e_ref{ifactorrep_e}.
Thus, the restriction map on the equivariant cohomology induced by the inclusion 
$P_i\!\lra\!\P^{n-1}$ is given~by
\begin{equation}\label{restrmap_e}
H_{\T}^*(\P^{n-1})=\Q[x,\al_1,\ldots,\al_n]\big/\prod_{k=1}^{k=n}(x\!-\!\al_k)
\lra H_{\T}^*(P_i)=\Q[\al_1,\ldots,\al_n], \qquad x\lra\al_i.
\end{equation}
By~\e_ref{restrmap_e},
\begin{equation}\label{uniquecond_e}
\eta=0 \in H_{\T}^*(\P^{n-1}) \qquad\Llra\qquad
\eta|_{P_i}=0\in H_{\T}^* ~~\forall~i=1,2,\ldots,n.
\end{equation}\\

\noindent
The tautological line bundle $\ga_{n-1}\!\lra\!\P^{n-1}$ is a subbundle of 
$\P^{n-1}\!\times\!\C^n$ preserved by the diagonal action of~$\T$.
Thus, the action of $\T$ on $\P^{n-1}$ naturally lifts to an action 
on~$\ga_{n-1}$ and
\begin{equation}\label{garestr_e}
\E\big(\ga_{n-1}^*\big)\big|_{P_i}=\al_i \qquad\forall~i=1,2,\ldots,n.
\end{equation}
Via the exact sequence
$$0\lra \ga_{n-1}^*\otimes\ga_{n-1}
\lra \ga_{n-1}^*\otimes\big(\P^{n-1}\!\times\!\C^n\big) 
\lra T\P^{n-1} \lra0$$
of vector bundles on $\P^{n-1}$, $\T$ also lifts to an action on $T\P^{n-1}$.
By~\e_ref{weightspord_e} and~\e_ref{garestr_e}, 
\begin{equation}\label{tangrestr_e}
\E\big(T\P^{n-1}\big)\big|_{P_i}=\prod_{k\neq i}(\al_i\!-\!\al_k)=\phi_i|_{P_i}
\qquad\forall~i=1,2,\ldots,n.
\end{equation}\\

\noindent
If $G$ acts smoothly on a smooth compact oriented manifold $M$, there is a well-defined
integration-along-the-fiber homomorphism
$$\int_M\!: H_G^*(M)\lra H_G^*$$
for the fiber bundle $BM\!\lra\!BG$.
The classical localization theorem of~\cite{ABo} relates it to 
integration along the fixed locus of the $G$-action.
The latter is a union of smooth compact orientable manifolds~$F$
and $G$ acts on the normal bundle $\N F$ of each~$F$.
Once an orientation of $F$ is chosen, there is a well-defined 
integration-along-the-fiber homomorphism
$$\int_F\!: H_G^*(F)\lra H_G^*.$$
The localization theorem states that 
\begin{equation}\label{ABothm_e}
\int_M\psi = \sum_F\int_F\frac{\psi|_F}{\E(\N F)} \in \H^*_G
\qquad\forall~\psi\in H_G^*(M),
\end{equation}
where the sum is taken over all components $F$ of the fixed locus of $G$. 
Part of the statement of~\e_ref{ABothm_e} is that $\E(\N F)$
is invertible in~$\H_G^*(F)$.
In the case of the standard action of $\T$ on~$\P^{n-1}$, \e_ref{ABothm_e} implies that 
\begin{equation}\label{phiprop_e}
\eta|_{P_i}=\int_{\P^{n-1}}\eta\phi_i \in\H_{\T}
\qquad\forall~\eta\!\in\!\H_{\T}^*(\P^{n-1}),~i=1,2,\ldots,n;
\end{equation}
see also \e_ref{tangrestr_e}.\\

\noindent
Finally, if $f\!:M\!\lra\!M'$ is a $G$-equivariant map between two compact oriented 
manifolds, there is a well-defined pushforward homomorphism
$$f_*\!: H_G^*(M) \lra H_G^*(M').$$
It is characterized by the property that 
\begin{equation}\label{pushdfn_e}
\int_{M'}(f_*\eta)\,\psi=\int_M\eta\,(f^*\psi) 
\qquad~\forall~\eta\in H_G^*(M), \psi\in H_G^*(M').
\end{equation}
The homomorphism $\int_M$ of the previous paragraph corresponds to $M'$ being a point.
It is immediate from~\e_ref{pushdfn_e} that 
\begin{equation}\label{pushprop_e}
f_*\big(\eta\,(f^*\psi)\big)=(f_*\eta)\,\psi 
\qquad~\forall~\eta\in H_G^*(M), \psi\in H_G^*(M').
\end{equation}

\subsection{Proof of Lemma~\ref{recgen_lmm}}
\label{local_subs1}

\noindent
The standard $\T$-action on $\P^{n-1}$ (as well as any other action) induces 
$\T$-actions on moduli spaces of stable maps $\ov\M_{0,m}(\P^{n-1},d)$ 
by composition on the right:
$$g\cdot[\cC,f]=[\cC,g\circ f]  \qquad\forall~g\in\T, 
\,[\cC,f]\in\ov\M_{0,m}(\P^{n-1},d).$$
All evaluation maps
$$\ev_i\!: \ov\M_{0,m}(\P^{n-1},d)\lra\P^{n-1}, \qquad
[\cC,y_1,\ldots,y_k,f]\lra f(y_i),$$
are $\T$-equivariant.
Via the natural lift of the $\T$-action to $\ga_{n-1}^*\!\lra\!\P^{n-1}$
described in Subsection~\ref{equivcoh_subs}, the $\T$-action on $\ov\M_{0,m}(\P^{n-1},d)$ 
lifts to $\T$-actions on the vector bundles~$\V_0$, $\V_0$, and~$\V_0''$,
as well as on the universal tangent line bundles.\\

\begin{figure}
\begin{pspicture}(-.7,-1.5)(10,2)
\psset{unit=.4cm}
\psline[linewidth=.04](2.5,0)(6,0)\rput(4.2,.7){\smsize{$d$}}
\pscircle*(2.5,0){.2}\rput(2.5,.7){\smsize{$i$}}
\pscircle*(6,0){.2}\rput(6,.75){\smsize{$j$}}
\psline[linewidth=.04](2.5,0)(1,1.5)\rput(1,2){\smsize{$\bf 1$}}
\psline[linewidth=.04](6,0)(7.5,1.5)\rput(8,1.5){\smsize{$\bf 2$}}
\psline[linewidth=.04](6,0)(7.5,-1.5)\rput(8,-1.5){\smsize{$\bf 3$}}
\rput(4.2,-1.5){$\Ga_{ij}(d)$}
\psline[linewidth=.04](11.5,0)(15,0)\rput(13.2,.7){\smsize{$d$}}
\pscircle*(11.5,0){.2}\rput(11.5,.7){\smsize{$i$}}
\pscircle*(15,0){.2}\rput(15,.75){\smsize{$j$}}
\psline[linewidth=.04](11.5,0)(10,1.5)\rput(10,2){\smsize{$\bf 1$}}
\psline[linewidth=.04](15,0)(17,2)\pscircle*(17,2){.2}
\psline[linewidth=.04](17,2)(20,3)\pscircle*(20,3){.2}
\psline[linewidth=.04](17,2)(20,1)\pscircle*(20,1){.2}
\psline[linewidth=.04](15,0)(17,-2)\pscircle*(17,-2){.2}
\psline[linewidth=.04](17,-2)(18.5,-.5)\rput(19,-.5){\smsize{$\bf 2$}}
\psline[linewidth=.04](17,-2)(18.5,-3.5)\rput(19,-3.5){\smsize{$\bf 3$}}
\psline[linewidth=.04](26.5,0)(30,0)\pscircle*(30,0){.2}
\pscircle*(26.5,0){.2}\rput(26.5,.7){\smsize{$i$}}
\psline[linewidth=.04](26.5,0)(25,1.5)\rput(25,2){\smsize{$\bf 1$}}
\psline[linewidth=.04](30,0)(32,2)\pscircle*(32,2){.2}
\psline[linewidth=.04](32,2)(35,3)\pscircle*(35,3){.2}
\psline[linewidth=.04](32,2)(35,1)\pscircle*(35,1){.2}
\psline[linewidth=.04](30,0)(32,-2)\pscircle*(32,-2){.2}
\psline[linewidth=.04](32,-2)(33.5,-.5)\rput(34,-.5){\smsize{$\bf 2$}}
\psline[linewidth=.04](26.5,0)(24.5,-2)\pscircle*(24.5,-2){.2}
\psline[linewidth=.04](26.5,0)(28.5,-2)\pscircle*(28.5,-2){.2}
\psline[linewidth=.04](26.5,0)(28.5,2)\pscircle*(28.5,2){.2}
\psline[linewidth=.04](28.5,2)(31.5,3)\pscircle*(31.5,3){.2}
\psline[linewidth=.04](32,-2)(33.5,-3.5)\rput(34,-3.5){\smsize{$\bf 3$}}
\end{pspicture}
\caption{Two graphs of type $A_i(j;d)\!\subset\!A_i$ and a graph of type $B_i$}
\label{loopgraph_fig}
\end{figure}
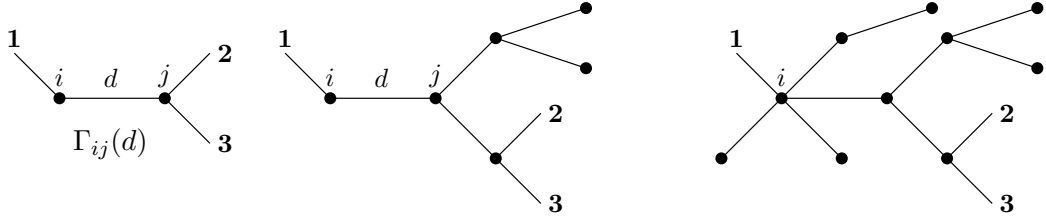

\noindent
As described in detail in \cite[Section~27.3]{MirSym},
the fixed loci of the $\T$-action on $\ov\M_{0,m}(\P^{n-1},d)$ 
are indexed by {\tt decorated graphs} that have no loops.
A {\tt graph} consists of a set~$\Ver$ of {\tt vertices} and 
a collection $\Edg$ of {\tt edges}, i.e.~of two-element subsets of~$\Ver$.
A {\tt loop} in a graph $(\Ver,\Edg)$ is a subset of $\Edg$ of the form
$$\big\{\{v_1,v_2\},\{v_2,v_3\},\ldots,\{v_N,v_1\}\big\},
\qquad v_1,\ldots,v_N\!\in\!\Ver,~ N\!\ge\!1.$$
Neither of the three graphs in Figure~\ref{loopgraph_fig} has a loop.
A {\tt decorated graph} is a tuple
\begin{equation}\label{decortgraphdfn_e}
\Ga = \big(\Ver,\Edg;\mu,\d,\eta\big),
\end{equation}
where $(\Ver,\Edg)$ is a graph and
$$\mu\!:\Ver\lra [n], \qquad \d\!: \Edg\lra\Z^+, \quad\hbox{and}\quad \eta\!:[m]\lra\Ver$$
are maps such that
\begin{equation}\label{decorgraphcond_e}
\mu(v_1)\neq\mu(v_2)  \qquad\hbox{if}\quad \{v_1,v_2\}\in\Edg.
\end{equation}
In Figure~\ref{loopgraph_fig}, the values of the map $\mu$ on some of the vertices
are indicated by letters next to those vertices.
Similarly, the value of the map $\d$ on one of the edges is indicated by a letter 
next to the edge.
The three elements of the~set $[m]\!=\![3]$ are shown in bold face.
They are linked by line segments to their images under~$\eta$.
By~\e_ref{decorgraphcond_e}, no two consecutive vertex labels are the same
and thus $j\!\neq\!i$.\\

\noindent
The fixed locus $\cZ_{\Ga}$ of $\ov\M_{0,m}(\Pn,d)$ corresponding to a 
decorated graph $\Ga$ consists of the stable maps~$f$ from a 
genus-zero nodal curve $\cC_f$ with $k$ marked points into~$\P^{n-1}$ 
that satisfy the following conditions.
The components of $\cC_f$ on which the map $f$ is not constant
are rational and correspond to the edges of~$\Ga$.
Furthermore, if $e\!=\!\{v_1,v_2\}$ is an edge,
the restriction of $f$ to the component $\cC_{f,e}$ corresponding to~$e$
is a degree-$\d(e)$ cover of the line 
$$\P^1_{\mu(v_1),\mu(v_2)}\subset\P^{n-1}$$
passing through the fixed points $P_{\mu(v_1)}$ and $P_{\mu(v_2)}$.
The map $u|_{\cC_{f,e}}$ is ramified only over $P_{\mu(v_1)}$ and~$P_{\mu(v_2)}$.
In particular, $f|_{\cC_{f,e}}$ is unique up to isomorphism.
The remaining, contracted, components of $\cC_f$ are  
indexed by the vertices $v\!\in\!\Ver$ such that 
$$\val(v)\equiv \big|\{v'\!\in\!\Ver\!: \{v,v'\}\!\in\!\Edg\}\big|
+ \big|\{i\!\in\![m]\!: \eta(i)\!=\!v\}\big| \ge 3.$$
The map $f$ takes such a component $\cC_{f,v}$ to the fixed point~$\mu(v)$.
Thus,
\begin{equation}\label{Zlocus_e}
\cZ_{\Ga}\approx \ov\cM_{\Ga}\!\equiv\!\prod_{v\in\Ver}\!\!\ov\cM_{0,\val(v)},
\end{equation}
where $\ov\cM_{0,l}$ denotes the moduli space
of stable genus-zero curves with $l$ marked points.
For the purposes of this definition, $\ov\cM_{0,1}$ and $\ov\cM_{0,2}$
are one-point spaces.
For example, in the case of the last diagram in Figure~\ref{loopgraph_fig},
$$\cZ_{\Ga}\approx
\ov\cM_{\Ga} \!\equiv \ov\cM_{0,5} \!\times\! \ov\cM_{0,3}^3
\!\times\! \ov\cM_{0,2} \!\times\! \ov\cM_{0,1}^5 \approx \ov\cM_{0,5}$$
is a fixed locus\footnote{after dividing by an appropriate automorphism 
group; see \cite[Section 27.3]{MirSym}}
in $\ov\M_{0,3}(\P^{n-1},d)$ for some $d\!\ge\!9$.\\

\noindent
We now verify Lemma~\ref{recgen_lmm}.
Let
\begin{equation}\label{etabedfn_e}
\eta^{\be}=\prod_{j=2}^{j=m}\big(\psi_j^{\be_j}\ev_j^*\eta_j\big).
\end{equation}
Suppose $\Ga$ is a decorated graph as in~\e_ref{decortgraphdfn_e} that contributes
to~\e_ref{recgen_e}, in the sense of the localization formula~\e_ref{ABothm_e}.
By~\e_ref{phidfn_e} and~\e_ref{restrmap_e},
$$\ev_1^*\phi_i=\prod_{k\neq i}\big(\al_{\mu(\eta(1))}-\al_k\big)
=\de_{i,\mu(\eta(1))}\prod_{k\neq i}(\al_i-\al_k),$$
where $\de_{i,\mu(\eta(1))}$ is the Kronecker delta function.
Thus, by~\e_ref{ABothm_e}, $\Ga$ does not contribute to~\e_ref{Z2wptdfn_e1}
unless $\mu(\eta(1))\!=\!i$, i.e.~the marked point of the map is taken 
to the point $P_i\!\in\!\P^{n-1}$.
There are two types of  graphs that do (or may) contribute to~\e_ref{Z2wptdfn_e1}; 
they will be called {\tt $A_i$} and {\tt $B_i$-types}.
In a graph of the $A_i$-type, the marked point~$1$ is attached to a vertex $v_0\!\in\!\Ver$ 
of valence two which is labeled~$i$.
In a graph of the $B_i$-type, the marked point~$1$ is attached to 
a vertex $v_0$ of valence at~$3$, which is still labeled~$i$. 
Examples of the two types are depicted in Figure~\ref{loopgraph_fig}.\\

\noindent
Suppose $\Ga$ is a graph of type $B_i$ and 
$$\cZ_{\Ga}\subset\ov\M_{0,m}(\P^{n-1},d),$$
so that $\Ga$ contributes to the coefficient of~$u^d$ in~\e_ref{Z2wptdfn_e1}.
In this case, the restriction of $\psi_1$ to $\cZ_{\Ga}$
is the pull-back of a $\psi$-class from the component $\ov\cM_{0,\val(v_0)}$
in the decomposition~\e_ref{Zlocus_e}.
Since the $\T$-action on the corresponding tautological line bundle is trivial,
$$\psi_1^k\big|_{\cZ_{\Ga}}=0 \qquad\forall~ k\ge d\!+\!m>\val(v_0)-3.$$
Thus, $\Ga$ contributes a polynomial in $\hb^{-1}$, of degree at most $d\!+\!m$,
to the coefficient of~$u^d$ in~\e_ref{Z2wptdfn_e1}.
Therefore, the contributions of the loci of type $B_i$ to~\e_ref{Z2wptdfn_e1}
are accounted for by the middle term in~\e_ref{recur_dfn_e}.\\

\noindent
A graph $\Ga$ as in~\e_ref{decortgraphdfn_e} of type $A_i$ has a unique vertex $v$
joined to~$v_0$. 
Denote by $A_i(j;d_0)$ the set of all graphs $\Ga$ of type~$A_i$ such that 
$\mu(v)\!=\!j$ and $\d(\{v_0,v\})\!=\!d_0$, i.e.~the unique vertex~$v$ of $\Ga$ 
joined to~$v_0$ is mapped to $P_j\!\in\!\P^{n-1}$ and the edge $\{v_0,v\}$ 
corresponds to the $d_0$-fold cover of $\P_{ij}^1$ branched only over $P_i$ and~$P_j$.
By~\e_ref{decorgraphcond_e},
\begin{equation}\label{Garestr_e}
A_i=\bigcup_{d_0=1}^{\i}\bigcup_{j\neq i}A_i(j;d_0).
\end{equation}\\

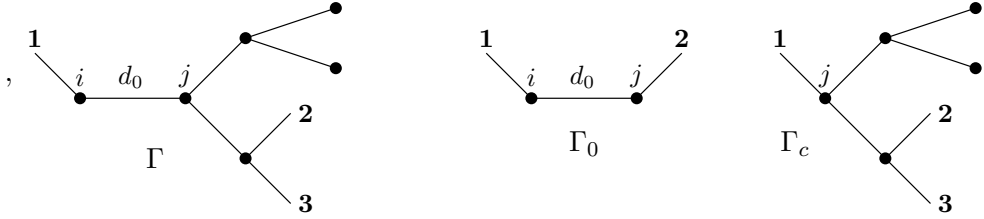
\begin{figure}
\begin{pspicture}(-1.2,-1.5)(10,2)
\psset{unit=.4cm}
\psline[linewidth=.04](2.5,0)(6,0)\rput(4.2,.7){\smsize{$d_0$}}
\pscircle*(2.5,0){.2}\rput(2.5,.7){\smsize{$i$}}
\pscircle*(6,0){.2}\rput(6,.75){\smsize{$j$}}
\psline[linewidth=.04](2.5,0)(1,1.5)\rput(1,2){\smsize{$\bf 1$}}
\psline[linewidth=.04](6,0)(8,2)\pscircle*(8,2){.2}
\psline[linewidth=.04](8,2)(11,3)\pscircle*(11,3){.2}
\psline[linewidth=.04](8,2)(11,1)\pscircle*(11,1){.2}
\psline[linewidth=.04](6,0)(8,-2)\pscircle*(8,-2){.2}
\psline[linewidth=.04](8,-2)(9.5,-.5)\rput(10,-.5){\smsize{$\bf 2$}}
\psline[linewidth=.04](8,-2)(9.5,-3.5)\rput(10,-3.5){\smsize{$\bf 3$}}
\rput(5,-2){$\Ga$}
\psline[linewidth=.04](17.5,0)(21,0)\rput(19.2,.7){\smsize{$d_0$}}
\pscircle*(17.5,0){.2}\rput(17.5,.7){\smsize{$i$}}
\pscircle*(21,0){.2}\rput(21,.75){\smsize{$j$}}
\psline[linewidth=.04](17.5,0)(16,1.5)\rput(16,2){\smsize{$\bf 1$}}
\psline[linewidth=.04](21,0)(22.5,1.5)\rput(22.5,2){\smsize{$\bf 2$}}'
\rput(19,-1.5){$\Ga_0$}
\pscircle*(27,0){.2}\rput(27,.75){\smsize{$j$}}
\psline[linewidth=.04](27,0)(25.5,1.5)\rput(25.5,2){\smsize{$\bf 1$}}
\psline[linewidth=.04](27,0)(29,2)\pscircle*(29,2){.2}
\psline[linewidth=.04](29,2)(32,3)\pscircle*(32,3){.2}
\psline[linewidth=.04](29,2)(32,1)\pscircle*(32,1){.2}
\psline[linewidth=.04](27,0)(29,-2)\pscircle*(29,-2){.2}
\psline[linewidth=.04](29,-2)(30.5,-.5)\rput(31,-.5){\smsize{$\bf 2$}}
\psline[linewidth=.04](29,-2)(30.5,-3.5)\rput(31,-3.5){\smsize{$\bf 3$}}
\rput(26,-1.5){$\Ga_c$}
\end{pspicture}
\caption{A graph of type $A_i^*(j;d_0)$ and its two subgraphs}
\label{splitgraph_fig}
\end{figure}

\noindent
Suppose $\Ga\!\in\!A_i(j;d_0)$ and $v$ is the unique vertex joined to $v_0$
by an edge.
We break at $v$ into two graphs:
\begin{enumerate}[label=(\roman*)]
\item $\Ga_0$ consisting of the vertices $v_0$ and $v$, the edge $\{v_0,v\}$, 
and marked points $1$ and $2$ attached to $v_0$ and $v$, respectively; 
\item $\Ga_c$ consisting all vertices and edges of $\Ga$, other than 
the vertex $v_0$ and the edge $\{v_0,v\}$, with a new marked point attached to~$v$;
\end{enumerate}
see Figure~\ref{splitgraph_fig}.
Let $d_c$ denote the degree of $\Ga_c$, i.e.~the sum of all edge labels.
By~\e_ref{Zlocus_e},
\begin{equation}\label{Zlocus_e2}
\cZ_{\Ga}\approx \cZ_{\Ga_0}\times\cZ_{\Ga_c}.
\end{equation}
Denote by $\pi_0$ and $\pi_c$ the two component projection maps.\\

\noindent
By Section 27.4 in~\cite{MirSym},
\begin{equation}\label{bndlsplit_e1}\begin{split}
\V_0''\big|_{\cZ_{\Ga}}&=\pi_0^*\V_0''\oplus\pi_c^*\V_0'',\\
\frac{\N\cZ_{\Ga}}{T_{P_i}\P^{n-1}}&=
\pi_0^*\bigg(\frac{\N\cZ_{\Ga_0}}{T_{P_i}\P^{n-1}}\bigg)
\oplus \pi_c^*\bigg(\frac{\N\cZ_{\Ga_c}}{T_{P_j}\P^{n-1}}\bigg)
\oplus \pi_0^*L_2\otimes\pi_c^*L_1,
\end{split}\end{equation}
where $L_2\!\lra\!\cZ_{\Ga_0}$ and $L_1\!\lra\!\cZ_{\Ga_c}$ are the tautological tangent
line bundles.
Thus, by~\e_ref{tangrestr_e},
\begin{equation}\label{bndlsplit_e2}\begin{split}
\frac{\E(\V_0'')\eta^{\be}}{\hb\!-\!\psi_1}\Big|_{\cZ_{\Ga}}
&=\pi_0^*\bigg(\frac{\E(\V_0'')}{\hb\!-\!\psi_1}\bigg)\cdot
\pi_c^*\Big(\E(\V_0'')\eta^{\be}\Big),\\
\frac{\ev_1^*\phi_i|_{\cZ_{\Ga}}}{\E(N\cZ_{\Ga})}&=
\pi_0^*\bigg(\frac{\ev_1^*\phi_i}{\E(\N\cZ_{\Ga_0})}\bigg)
\cdot \pi_c^*\bigg(\frac{\ev_1^*\phi_j}{\E(\N\cZ_{\Ga_0})}\bigg)
\cdot\frac{1}{\pi_0^*c_1(L_2)-\pi_c^*\psi_1}.
\end{split}\end{equation}
By Sections 27.1 and 27.2 in~\cite{MirSym}, on $\cZ_{\Ga_0}$
\begin{equation}\label{bndlsplit_e3}\begin{split}
\E(\V_0'')&=\prod_{r=0}^{nd_0-1}
\Big(n\al_i+r\frac{\al_j\!-\!\al_i}{d_0}\Big), \qquad
\psi_1=c_1(L_2)=\frac{\al_j\!-\!\al_i}{d_0},\\
\E(N\cZ_{\Ga_0})&=(-1)^{d_0}\prod_{r=1}^{r=d_0}\Big(r\frac{\al_j\!-\!\al_i}{d_0}\Big)^2
\prod_{r=0}^{r=d_0}\prod_{k\neq i,j}
\Big(\al_i\!-\!\al_k\!+\!r\frac{\al_j\!-\!\al_i}{d_0}\Big).
\end{split}\end{equation}
Thus, using~\e_ref{tangrestr_e} and taking into the account the automorphism group, 
$\Z_{d_0}$, we obtain
\begin{equation}\label{bndlsplit_e5}
\int_{\cZ_{\Ga_0}}\frac{\E(\V_0'')\ev_1^*\phi_i}{(\hb\!-\!\psi_1)\E(N\cZ_{\Ga_0})}
=\frac{1}{\hb-\frac{\al_j-\al_i}{d_0}}C_i^j(d_0).
\end{equation}
By \e_ref{Zlocus_e2}, \e_ref{bndlsplit_e3}, and~\e_ref{bndlsplit_e5}, 
the contribution of $\Ga$ to~\e_ref{recgen_e} is 
\begin{equation}\label{bndlsplit_e6}\begin{split}
&u^{d_0+d_c}\!\int_{\cZ_{\Ga}}
\frac{\E(\V_0'')\ev_1^*\phi_i\eta^{\be}}{\hb\!-\!\psi_1}\Big|_{\cZ_{\Ga}}
\frac{1}{\E(N\cZ_{\Ga})}\\
&\hspace{1.2in}
=\frac{C_i^j(d_0)u^{d_0}}{\hb-\frac{\al_j-\al_i}{d_0}}\cdot
\Bigg(\bigg\{u^{d_c}\!\int_{\cZ_{\Ga}}
\frac{\E(\V_0'')\ev_1^*\phi_j\eta^{\be}}{\hb\!-\!\psi_1}
\frac{1}{\E(N\cZ_{\Ga_c})}\bigg\}\bigg|_{\hb=\frac{\al_j-\al_i}{d_0}}\Bigg).
\end{split}\end{equation}\\

\noindent
We next sum~\e_ref{bndlsplit_e6} over $\Ga\!\in\!A_i(j;d_0)$.
This  is the same as summing the expression in the curly brackets over
all $m$-pointed graphs with the marked point~$1$ attached to a vertex~$v$ labeled~$j$,
i.e.~all graphs of types~$A_j$ and~$B_j$.
By the localization formula~\e_ref{ABothm_e}, the sum of the terms in the curly brackets
over all such graphs~$\Ga_c$ is $\cZ_{\eta,\be}(\hb,\al_j,u)$.
Thus,
\begin{equation}\label{bndlsplit_e7}\begin{split}
&\sum_{\Ga\in A_i(j;d_0)}u^{d_0+d_c}\!\!
\int_{\cZ_{\Ga}}\frac{\E(\V_0'')\ev_1^*\phi_i\eta^{\be}}{\hb\!-\!\psi_1}\Big|_{\cZ_{\Ga}}
\frac{1}{\E(N\cZ_{\Ga})}
=\frac{C_i^j(d_0)u^{d_0}}{\hb-\frac{\al_j-\al_i}{d_0}}\cdot
\cZ_{\eta,\be}\big((\al_j-\al_i)/d_0,\al_j,u\big).
\end{split}\end{equation}
We conclude that $\cZ_{\eta,\be}(\hb,x,u)$ is $C$-recursive in 
the sense of Definition~\ref{recur_dfn}:
\begin{itemize}
\item the middle term in~\e_ref{recur_dfn_e} consists of 
the contributions from the graphs of type~$B_i$;
\item the $(d_0,j)$-summand in~\e_ref{recur_dfn_e} consists of
the contributions from the graphs of type~$A_i(j;d_0)$.
\end{itemize}

\subsection{Proof of Lemma~\ref{polgen_lmm}}
\label{local_subs2}

\noindent
In this subsection we deduce Lemma~\ref{polgen_lmm} from Lemma~\ref{PhiZstr_lmm},
which is proved in the next subsection.
The argument, in this subsection and the next one, is a modification 
on the proof of self-polynomiality of~$\cZ$ in Section~30.2 of~\cite{MirSym}.\\

\noindent
We will denote the weight of the standard action of the one-torus $\T^1$ on $\C$ by~$\hb$.
Thus, by Subsection~\ref{equivcoh_subs},
$$H_{\T^1}^*\approx\Q[\hb], \quad H_{\wt\T}^*\approx\Q[\hb,\al_1,\ldots,\al_n]
\qquad\Lra\qquad \H_{\wt\T}^*\approx\Q_{\al}(\hb).$$
Throughout this subsection, $V\!=\!\C\!\oplus\!\C$ will denote the representation
of $\T^1$ with the weights $0$ and~$-\hb$.
The induced action on $\P V$ has two fixed points:
$$q_1\equiv[1,0], \qquad q_2\equiv[0,1].$$
Let $\ga_1\!\lra\!\P V$ be the tautological line bundle.
Then,
\begin{equation}\label{hbaract_e}
\E(\ga_1^*)\big|_{q_1}=0, \quad \E(\ga_1^*)\big|_{q_2}=-\hb, 
\quad \E(T_{q_1}\P V)=\hb, \quad \E(T_{q_2}\P V)=-\hb.
\end{equation}\\

\noindent
For each $d\!\in\!\bar\Z^+$, the action of $\wt\T$ on $\C^n\!\otimes\!\Sym^dV^*$ induces 
an action on 
$$\ov\X_d\equiv\P\big(\C^n\!\otimes\!\Sym^dV^*\big).$$
It has $(d\!+\!1)n$ fixed points:
$$P_i(r)\equiv \big[\ti{P}_i\otimes u^{d-r}v^r\big], \qquad
i\in[n],~r\in\{0\}\!\cup\![d],$$
if $(u,v)$ are the standard coordinates on $V$ and $\ti{P}_i\!\in\!\C^n$ is
the $i$th coordinate vector (so that $[\ti{P}_i]\!=\!P_i\!\in\!\P^{n-1}$).
Let
$$\Om\equiv \E(\ga^*)\in H_{\wt\T}^*\big(\ov\X_d\big)$$
denote the equivariant hyperplane class.\\

\noindent
For all $i\!\in\![n]$ and $r\in\{0\}\!\cup\![d]$, 
\begin{equation}\label{Xrestr_e}
\Om|_{P_i(r)}=\al_i\!+\!r\hb, \qquad 
\E(T_{P_i(r)}\ov\X_d)=\Bigg\{\underset{(s,k)\neq(r,i)}{\prod_{s=0}^{s=d}\prod_{k=1}^{k=n}}
(\Om\!-\!\al_k\!-\!r\hb)\Bigg\}\bigg|_{\Om=\al_i+r\hb}.\footnotemark
\end{equation}
\footnotetext{The weight (i.e.~negative first chern class) of the $\ti{T}$-action 
on the line $P_i(r)\!\subset\!\C^n\!\otimes\!\Sym^dV^*$ is $\al_i\!+\!r\hb$.
The tangent bundle of $\ov\X_d$ at $P_i(r)$ is the direct sum of the lines
$P_i(r)^*\!\otimes\!P_k(s)$ with $(k,s)\!\neq\!(i,r)$.}
Since 
\begin{gather*}
B\ov\X_d=\P\big(B(\C^n\!\otimes\!\Sym^dV^*)\big)\lra B\wt\T \qquad\hbox{and}\\
c\big(B(\C^n\!\otimes\!\Sym^dV^*)\big)
=\prod_{s=0}^{s=d}\prod_{k=1}^{k=n}\big(1-(\al_k\!+\!s\hb)\big)
\in H^*\big(B\wt\T),\footnotemark
\end{gather*}
\footnotetext{The vector space $\C^n\!\otimes\!\Sym^dV^*$ is the direct sum of
the one-dimensional representations $P_k(s)$ of~$\wt\T$.}
the $\wt\T$-equivariant cohomology of $\ov\X_d$ is given~by
\begin{equation}\label{Xdcohom_e}\begin{split}
H_{\ti\T}^*\big(\ov\X_d\big)&\equiv H^*\big(B\ov\X_d\big)
=H^*\big(B\wt\T\big)\big[\Om\big]\Big/
 \prod_{s=0}^{s=d}\prod_{k=1}^{k=n}\big(\Om-(\al_k\!+\!s\hb)\big)\\
&\approx \Q\big[\Om,\hb,\al_1,\ldots,\al\big]\Big/
 \prod_{s=0}^{s=d}\prod_{k=1}^{k=n}\big(\Om-\al_k-\!s\hb\big)\\
&\subset \Q_{\al}[\hb,\Om]\Big/
 \prod_{s=0}^{s=d}\prod_{k=1}^{k=n}\big(\Om-\al_s-r\hb\big).
\end{split}\end{equation}\\

\noindent
There is a natural $\wt\T$-equivariant morphism
$$\Th\!:\ov\M_{0,m}\big(\P V\!\times\!\P^{n-1},(1,d)\big)\lra \ov\X_d.$$
A general element of $b$ of $\ov\M_{0,m}\big(\P V\!\times\!\P^{n-1},(1,d)\big)$
determines a map
$$(f,g)\!:\P^1\lra(\P V,\P^n),$$ 
up to an automorphism  of the domain~$\P^1$.
Thus, the map
$$g\circ f^{-1}\!: \P V\lra \P^{n-1}$$
is well-defined and determines an element $\Th(b)\!\in\!\ov\X_d$.
The map $\Th$ extends continuously over the boundary of 
$\ov\M_{0,m}\big(\P V\!\times\!\P^{n-1},(1,d)\big)$.\footnote{For a complete algebraic proof,
see Lemma~2.6 in~\cite{LLY}.}
We denote the restriction of~$\Th$ to the smooth substack
\begin{equation}\label{Xdfn_e}
\X_d\equiv\big\{b\!\in\!\ov\M_{0,m}\big(\P V\!\times\!\P^{n-1},(1,d)\big)\!:
\ev_1(b)\!\in\!q_1\!\times\!\P^{n-1},~\ev_2(b)\!\in\!q_2\!\times\!\P^{n-1}\big\}
\end{equation}
of $\ov\M_{0,m}\big(\P V\!\times\!\P^{n-1},(1,d))$ by $\th_d$, or simply by $\th$
whenever there is no ambiguity.\\

\noindent
Let 
$$\pi\!:\ov\M_{0,m}\big(\P V\!\times\!\P^{n-1},(1,d)\big)\lra 
\ov\M_{0,m}\big(\P^{n-1},d\big)$$
be the natural projection map.

\begin{lmm}
\label{PhiZstr_lmm}
With $\cZ_{\eta,\be}$ as in Lemma~\ref{polgen_lmm} and $\Phi$ as in \e_ref{Phidfn_e},
\begin{equation}\label{PhiZstr_e}
(-\hb)^{m-2}\Phi_{\cZ,\cZ_{\eta,\be}}(\hb,u,z)=\sum_{d=0}^{\i}u^d\!\!
\int_{\X_d}\!\!\!e^{(\th^*\Om)z}
\pi^*\bigg(\E(\V_0'')\prod_{j=2}^{j=m}\big(\psi_j^{\be_j}\ev_j^*\eta_j\big)\bigg)
\prod_{j=3}^{j=m}\!\!\ev_j^*\big(\E(\ga_1^*)\big).
\end{equation}\\
\end{lmm}

\noindent
Similarly Section 30.2 in~\cite{MirSym}, this lemma implies that
$$(-\hb)^{m-2}\Phi_{\cZ,\cZ_{\eta,\be}}(\hb,u,z)
\in\Q_{\al}[\hb]\big[\big[u,z\big]\big]$$
for the following reason.
With $\eta^{\be}$ as in~\e_ref{etabedfn_e}, by~\e_ref{Xdcohom_e}
$$\th_{d*}\bigg(\pi^*\big(\E(\V_0'')\eta^{\be}\big)
\prod_{j=3}^{j=m}\ev_j^*\big(\E(\ga_1^*)\big)\bigg)
=E_{\cZ,\cZ_{\eta,\be};d}(\hb,\Om)$$
for some $E_{\cZ,\cZ_{\be,\eta};d}\!\in\!\Q_{\al}[\hb,\Om]$ 
of $\Om$-degree at most $(d\!+\!1)n\!-\!1$.
Therefore, by Lemma~\ref{PhiZstr_lmm}, \e_ref{pushprop_e}, \e_ref{ABothm_e},
and~\e_ref{Xrestr_e},
\begin{equation*}\begin{split}
(-\hb)^{m-2}\Phi_{\cZ,\cZ_{\eta,\be}}(\hb,u,z) &
=\sum_{d=0}^{\i}u^d\!\!\int_{\ov\X_d}e^{\Om z}
\th_{d*}\bigg(\pi^*\big(\E(\V_0'')\eta^{\be}\big)
\prod_{j=3}^{j=m}\ev_j^*\big(\E(\ga_1^*)\big)\bigg)\\
&=\sum_{d=0}^{\i}\!u^d\!\! \int_{\ov\X_d}\! e^{\Om z}E_{\cZ,\cZ_{\be,\eta};d}(\hb,\Om)
=\sum_{d=0}^{\i}\!u^d\Bigg(\sum_{r=0}^{r=d}\sum_{i=1}^{i=n}
\frac{e^{\Om z}E_{\cZ,\cZ_{\be,\eta};d}(\hb,\Om)|_{P_i(r)}}{\E(T_{P_i(r)}\ov\X_d)}\Bigg)\\
&=\sum_{d=0}^{\i}u^d\Bigg(\sum_{r=0}^{r=d}\sum_{i=1}^{i=n}
\frac{e^{\Om z}E_{\cZ,\cZ_{\be,\eta};d}(\hb,\Om)}
{\underset{(s,k)\neq(r,i)}{\prod_{s=0}^{s=d}\prod_{k=1}^{k=n}}(\Om\!-\!\al_k\!-\!s\hb)}
\bigg|_{\Om=\al_i+r\hb}\Bigg)\\
&=\sum_{d=0}^{\i}u^d\Bigg(\frac{1}{2\pi\I}\oint e^{\Om z}\frac{E_{\cZ,\cZ_{\be,\eta};d}(\hb,\Om)}
{\prod_{s=0}^{s=r}\prod_{k=1}^{k=n}(\Om\!-\!\al_k\!-\!s\hb)}d\Om\Bigg).
\end{split}\end{equation*}
In the last expression, the integral has the same meaning as in Lemma~\ref{Phistr_lmm1}.
We have thus shown that $\cZ_{\eta,\be}$ is polynomial with respect to~$\cZ$,
assuming Lemma~\ref{PhiZstr_lmm}.

\subsection{Proof of Lemma~\ref{PhiZstr_lmm}}
\label{PhiZstrpf_subs}

\noindent
In this subsection we use the localization formula~\e_ref{ABothm_e} to prove 
Lemma~\ref{PhiZstr_lmm}.
We show that each fixed locus of the $\wt\T$-action on $\X_d$ contributing
to the right-hand side of~\e_ref{PhiZstr_e} corresponds to a pair $(\Ga_1,\Ga_2)$ 
of a graphs, with $\Ga_1$ and $\Ga_2$ contributing to $\cZ(\hb,\al_i,u)$ and
$(-\hb)^{m-2}\cZ_{\eta,\be}(\hb,\al_i,u)$, respectively, for some $i\!\in\![n]$.\\

\noindent
Similarly to Subsection~\ref{local_subs1}, the fixed loci of the $\wt\T$-action on 
$\ov\M_{0,m}\big(\P V\!\times\!\P^{n-1},(d',d)\big)$
correspond to decorated graphs $\Ga$ with $m$ marked points and no loops.
Each edge should be labeled by a pair of integers, indicating the degrees of
the corresponding maps in $\P V$ and in~$\P^{n-1}$.
Each vertex should be labeled either $(1,j)$ or $(2,j)$ for some $j\!\in\![n]$,
indicating the fixed point, $(q_1,P_j)$ or $(q_2,P_j)$, to which the vertex is mapped.
No two consecutive vertex labels are the same, but if 
two consecutive vertex labels differ in the $k$th component, with $k\!=\!1,2$,
the $k$th component of the label for edge connecting them must be nonzero.\\

\begin{figure}
\begin{pspicture}(-.5,-1.2)(10,1.8)
\psset{unit=.4cm}
\psline[linewidth=.1](17,0)(22,0)
\pscircle*(17,0){.2}\rput(17,.7){\smsize{$i$}}
\pscircle*(22,0){.2}\rput(22,.7){\smsize{$i$}}
\psline[linewidth=.04](17,0)(14,1)\pscircle*(14,1){.2}
\rput(15.5,1){\smsize{$2$}}\rput(14,1.7){\smsize{$1$}}
\psline[linewidth=.04](14,1)(11.5,2.5)\pscircle*(11.5,2.5){.2}
\rput(12.8,2.3){\smsize{$2$}}\rput(10.9,2.5){\smsize{$2$}}
\psline[linewidth=.04](14,1)(11.5,-.5)\pscircle*(11.5,-.5){.2}
\rput(12.3,.6){\smsize{$1$}}\rput(10.9,-.5){\smsize{$3$}}
\psline[linewidth=.04](17,0)(14.5,-1.5)\pscircle*(14.5,-1.5){.2}
\rput(16.2,-1.1){\smsize{$3$}}\rput(14.5,-.8){\smsize{$3$}}
\psline[linewidth=.04](14.5,-1.5)(12,-2)\rput(11.5,-2){\smsize{$\bf 1$}}
\psline[linewidth=.04](22,0)(25,1)\pscircle*(25,1){.2}
\rput(23.5,1){\smsize{$1$}}\rput(25,1.7){\smsize{$4$}}
\psline[linewidth=.04](25,1)(27,2)\rput(27.5,2){\smsize{$\bf 2$}}
\psline[linewidth=.04](25,1)(27.5,-.5)\pscircle*(27.5,-.5){.2}
\rput(26.7,.6){\smsize{$3$}}\rput(28.1,-.5){\smsize{$1$}}
\psline[linewidth=.04](22,0)(24.5,-1.5)\pscircle*(24.5,-1.5){.2}
\rput(22.9,-1.1){\smsize{$5$}}\rput(24.5,-.8){\smsize{$1$}}
\psline[linewidth=.04](24.5,-1.5)(27,-2)\rput(27.5,-2){\smsize{$\bf 3$}}
\end{pspicture}
\caption{A graph representing a fixed locus in $\X_d$; $i\!\neq\!1,3,4$}
\label{Xlocus_fig}
\end{figure}
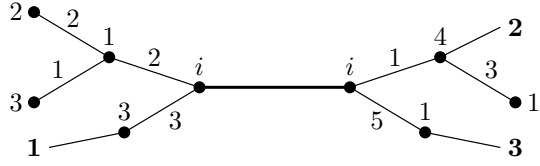

\noindent
The situation for the  $\wt\T$-action on 
$$\X_d\subset\ov\M_{0,m}\big(\P V\!\times\!\P^{n-1},(1,d)\big)$$
is simpler, however.
There is a unique edge of positive $\P V$-degree. 
We draw it as a thick horizontal line.
The first component of all other edge labels must be~$0$; so we drop~it.
The first components of the vertex labels change only when the thick edge is crossed.
Thus, we drop the first components of the vertex labels as well,
with the convention that these components are $1$ on the left side of the thick
edge and $2$ on the right.
In particular, the vertices to the left of the thick edge (including the left endpoint)
lie in $q_1\!\times\!\P^{n-1}$ and the vertices to its right lie in $q_2\!\times\!\P^{n-1}$.
Thus, by~\e_ref{Xdfn_e}, the marked point~$1$ is attached to a vertex to the left of
the thick edge and the marked point~$2$ is attached to a vertex to the right.
Furthermore, by the first identity in~\e_ref{hbaract_e}, such a graph will not
contribute to the right-hand side of~\e_ref{PhiZstr_e} unless 
the remaining marked points are also attached to vertices to the right of the thick edge. 
Finally, both  vertices of the thick edge have the same (remaining, second) label $i\!\in\![n]$.
Let $\cA_i$ denote the set of graphs as above so that the two endpoints of the thick
edge are labeled~$i$; see Figure~\ref{Xlocus_fig}.\\

\begin{figure}
\begin{pspicture}(-.5,-1.2)(10,1.8)
\psset{unit=.4cm}
\pscircle*(10,0){.2}\rput(10,.7){\smsize{$i$}}
\psline[linewidth=.04](10,0)(7,1)\pscircle*(7,1){.2}
\rput(8.5,1){\smsize{$2$}}\rput(7,1.7){\smsize{$1$}}
\psline[linewidth=.04](7,1)(4.5,2.5)\pscircle*(4.5,2.5){.2}
\rput(7.8,2.3){\smsize{$2$}}\rput(3.9,2.5){\smsize{$2$}}
\psline[linewidth=.04](7,1)(4.5,-.5)\pscircle*(4.5,-.5){.2}
\rput(5.3,.6){\smsize{$1$}}\rput(3.9,-.5){\smsize{$3$}}
\psline[linewidth=.04](10,0)(7.5,-1.5)\pscircle*(7.5,-1.5){.2}
\rput(9.2,-1.1){\smsize{$3$}}\rput(7.5,-.8){\smsize{$3$}}
\psline[linewidth=.04](7.5,-1.5)(5,-2)\rput(4.5,-2){\smsize{$\bf 1$}}
\psline[linewidth=.04](10,0)(11.5,1.5)\rput(11.5,2){\smsize{$\bf 2$}}
\psline[linewidth=.04](17,0)(15.5,1.5)\rput(15.5,2){\smsize{$\bf 1$}}
\psline[linewidth=.04](22,0)(23.5,1.5)\rput(23.5,2){\smsize{$\bf 2$}}
\psline[linewidth=.1](17,0)(22,0)
\pscircle*(17,0){.2}\rput(17,.7){\smsize{$i$}}
\pscircle*(22,0){.2}\rput(22,.7){\smsize{$i$}}
\psline[linewidth=.04](29,0)(27.5,1.5)\rput(27.5,2){\smsize{$\bf 1$}}
\pscircle*(29,0){.2}\rput(29,.7){\smsize{$i$}}
\psline[linewidth=.04](29,0)(32,1)\pscircle*(32,1){.2}
\rput(30.5,1){\smsize{$1$}}\rput(32,1.7){\smsize{$4$}}
\psline[linewidth=.04](32,1)(34,2)\rput(34.5,2){\smsize{$\bf 2$}}
\psline[linewidth=.04](32,1)(34.5,-.5)\pscircle*(34.5,-.5){.2}
\rput(33.7,.6){\smsize{$3$}}\rput(35.1,-.5){\smsize{$1$}}
\psline[linewidth=.04](29,0)(31.5,-1.5)\pscircle*(31.5,-1.5){.2}
\rput(29.9,-1.1){\smsize{$5$}}\rput(31.5,-.8){\smsize{$1$}}
\psline[linewidth=.04](31.5,-1.5)(34,-2)\rput(34.5,-2){\smsize{$\bf 3$}}
\end{pspicture}
\caption{The three sub-graphs of the graph in Figure~\ref{Xlocus_fig}}
\label{Xsplit_fig}
\end{figure}
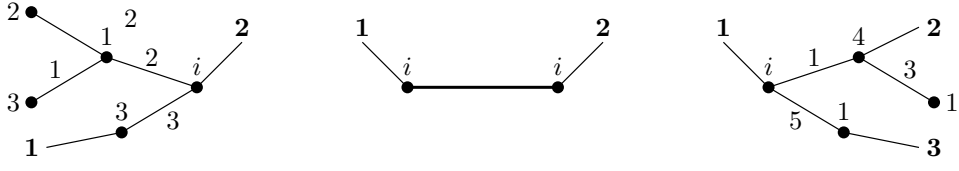

\noindent
If $\Ga\!\in\!\cA_i$, we break it into three sub-graphs:
\begin{enumerate}[label=(\roman*)]
\item $\Ga_1$ consisting of all vertices and edges of $\Ga$ to the left of the thick edge,
including its left vertex~$v_1$, and a new marked point~$2$ attached to~$v_1$;
\item $\Ga_0$ consisting of the thick edge, its two vertices $v_1$ and $v_2$, and 
new marked points $1$ and $2$ attached to $v_1$ and $v_2$, respectively;
\item $\Ga_2$ consisting of all vertices and edges of $\Ga$ to the right of the thick edge,
including its right vertex~$v_2$, and a new marked point~$1$ attached to~$v_2$;
\end{enumerate}
see Figure~\ref{Xsplit_fig}. The fixed locus in $\X_d$ corresponding to $\Ga$ is then  
\begin{equation}\label{Xlocus_e1}
\cZ_{\Ga}\approx \cZ_{\Ga_1}\times\cZ_{\Ga_0}\times\cZ_{\Ga_2}.
\end{equation}
The middle term is a single point.
Let $\pi_1$, $\pi_0$, and $\pi_2$ denote the three component projection maps.
Denote by $d_1$ and $d_2$ be the degrees of $\Ga_1$ and $\Ga_2$, i.e.
$$\cZ_{\Ga_1}\subset\ov\M_{0,2}(\P^{n-1},d_1),\qquad
\cZ_{\Ga_2}\subset\ov\M_{0,m}(\P^{n-1},d_2).$$
The exceptional case for the first statement is $d_1\!=\!0$, in which case
the corresponding moduli space does not exist.\\

\noindent
Suppose $\Ga\!\in\!\cA_i$, $d_1$ and $d_2$ are as above, and $d_1\!>\!0$.
Similarly to~\e_ref{bndlsplit_e1},
\begin{equation}\label{Xbndlsplit_e1}\begin{split}
\pi^*\V_0''\big|_{\cZ_{\Ga}}&=\pi_1^*\V_0''\oplus\pi_2^*\V_0'',\\
\frac{\N\cZ_{\Ga}}{T_{P_i}\P^{n-1}}&=
\pi_1^*\bigg(\frac{\N\cZ_{\Ga_1}}{T_{P_i}\P^{n-1}}\bigg)
\oplus \pi_2^*\bigg(\frac{\N\cZ_{\Ga_2}}{T_{P_i}\P^{n-1}}\bigg)
\oplus \pi_1^*L_2\otimes\pi_0^*L_1 \oplus \pi_0^*L_2\otimes\pi_2^*L_1,
\end{split}\end{equation}
where $\N\cZ_{\Ga}\!\lra\!\cZ_{\Ga}$ is the normal bundle of $\cZ_{\Ga}$ in $\X_d$
and $L_2\!\lra\!\cZ_{\Ga_1}$, $L_1,L_2\!\lra\!\cZ_{\Ga_0}$, and
$L_1\!\lra\!\cZ_{\Ga_2}$ are the tautological tangent line bundles.
We note that
$$L_1=T_{q_1}\P^1 \quad\hbox{and}\quad L_2=T_{q_2}\P^1 \qquad\hbox{on}\quad \cZ_{\Ga_0}.$$
Thus, by~\e_ref{Xbndlsplit_e1}, \e_ref{tangrestr_e}, and~\e_ref{hbaract_e},
\begin{alignat}{1}\label{Xbndlsplit_e3}
\pi^*\big(\E(\V_0'')\eta^{\be}\big)\prod_{j=3}^{j=m}\ev_j^*(\E(\ga^*))\bigg|_{\cZ_{\Ga}}
&=\pi_1^*\E(\V_0'')\cdot\pi_2^*\big(\E(\V_0'')\eta^{\be}(-\hb)^{m-2}\big),\\
\frac{\prod_{k\neq i}(\al_i\!-\!\al_k)}{\E(N\cZ_{\Ga})}&=
\pi_1^*\bigg(\frac{\ev_2^*\phi_i}{\E(\N\cZ_{\Ga_0})}\bigg)
\cdot \pi_2^*\bigg(\frac{\ev_1^*\phi_i}{\E(\N\cZ_{\Ga_0})}\bigg)
\cdot\frac{1}{\hb-\pi_1^*\psi_2}\cdot\frac{1}{(-\hb)-\pi_2^*\psi_1}.\notag
\end{alignat}
The map $\th$ takes the locus $\cZ_{\Ga}$ to a fixed point $P_k(r)\!\in\!\ov\X_d$.
It is immediate that $k\!=\!i$.
By continuity considerations, $r\!=\!d_1$.
Thus, by the first identity in \e_ref{Xrestr_e},
$$\th^*\Om\big|_{\cZ_{\Ga}}=\al_i\!+\!d_1\hb.$$
Combining \e_ref{Xlocus_e1} and~\e_ref{Xbndlsplit_e3} with this observation, we obtain
\begin{equation}\label{Xbndlsplit_e5}\begin{split}
&\int_{\cZ_{\Ga}}
\frac{e^{(\th^*\Om)z}\pi^*\big(\E(\V_0'')\eta^{\be}\big)\prod_{j=3}^{j=m}
\ev_j^*(\E(\ga^*)\big)|_{\cZ_{\Ga}}}{\E(\N\cZ_{\Ga})}
=(-\hb)^{m-2}\frac{e^{\al_iz}}{\prod_{k\neq i}(\al_i\!-\!\al_k)}\\
&\hspace{1in}
\times \Bigg\{e^{d_1\hb z}\!\! \int_{\cZ_{\Ga_1}} \!\!\!\!\!\!
\frac{\E(\V_0'')\ev_2^*\phi_i}{\hb\!-\!\psi_2}
\Big|_{\cZ_{\Ga_1}}\frac{1}{\E(\N\cZ_{\Ga_1})}\Bigg\}
\Bigg\{\int_{\cZ_{\Ga_2}} \!
\frac{\E(\V_0'')\eta^{\be}}{(-\hb)\!-\!\psi_1}
\Big|_{\cZ_{\Ga_2}}\frac{1}{\E(\N\cZ_{\Ga_2})}\Bigg\}.
\end{split}\end{equation}
We note that this identity remains valid for $d_1\!=\!0$ if we set 
the term in the first curly brackets to~$1$ for $d_1\!=\!0$.\\

\noindent
We now sum up \e_ref{Xbndlsplit_e5}, multiplied by $u^{d_1+d_2}$, over all 
$\Ga\!\in\!\cA_i$. 
This is the same as summing over all pairs $(\Ga_1,\Ga_2)$ of graphs such that 
\begin{enumerate}[label=(\arabic*)]
\item $\Ga_1$ is a 2-pointed graph of a degree $d_1\!\ge\!0$ such that the marked point $2$
is attached to the vertex labeled~$i$;
\item $\Ga_2$ is an $m$-pointed graph of a degree $d_2\!\ge\!0$ such that the marked point $1$
is attached to the vertex labeled~$i$.
\end{enumerate}
By the localization formula~\e_ref{ABothm_e},
\begin{equation}\label{Xbndlsplit_e7}\begin{split}
&\sum_{\Ga_1}u^{d_1}\Bigg\{e^{d_1\hb z}\!\! \int_{\cZ_{\Ga_1}} \!\!\!\!\!\!
\frac{\E(\V_0'')\ev_2^*\phi_i}{\hb\!-\!\psi_2}
\Big|_{\cZ_{\Ga_1}}\frac{1}{\E(\N\cZ_{\Ga_1})}\Bigg\}
=1+\sum_{d=1}^{\i}(ue^{\hb z})^d
\int_{\ov\M_{0,2}(\P^{n-1},d)}\frac{\E(\V_0'')}{\hb\!-\!\psi_2}\ev_2^*\phi_i\\
&\hspace{1in}
=1+\sum_{d=1}^{\i}(ue^{\hb z})^d
\int_{\ov\M_{0,2}(\P^{n-1},d)}\frac{\E(\V_0')}{\hb\!-\!\psi_1}\ev_1^*\phi_i
=\cZ\big(\hb,\al_i,ue^{\hb z}\big);\\
&\sum_{\Ga_2}u^{d_2}\Bigg\{\int_{\cZ_{\Ga_2}} \!\!\!\!\!\!
\frac{\E(\V_0'')\eta^{\be}\ev_1^*\phi_i}{(-\hb)\!-\!\psi_1}
\Big|_{\cZ_{\Ga_2}}\frac{1}{\E(\N\cZ_{\Ga_2})}\Bigg\}
=\sum_{d=0}^{\i}u^d
\int_{\ov\M_{0,m}(\P^{n-1},d)}\!\!\!
\frac{\E(\V_0'')\eta^{\be}\ev_1^*\phi_i}{(-\hb)\!-\!\psi_1}\\
&\hspace{3.55in} = \cZ_{\eta,\be}\big(-\hb,\al_i,u\big).
\end{split}\end{equation}
Finally, by~\e_ref{ABothm_e}, \e_ref{Xbndlsplit_e5}, and~\e_ref{Xbndlsplit_e7},
\begin{equation*}\begin{split}
\int_{\X_d}e^{(\th^*\Om)z}\pi^*\big(\E(\V_0'')\eta^{\be}\big)
\prod_{j=3}^{j=m}\ev_j^*\big(\E(\ga_1^*)\big)
&=(-\hb)^{m-2}\sum_{i=1}^{i=n}\frac{e^{\al_iz}}{\prod\limits_{k\neq i}(\al_i\!-\!\al_k)}
\cZ\big(\hb,\al_i,ue^{\hb z}\big)\cZ_{\eta,\be}\big(-\hb,\al_i,u\big)\\
&=(-h)^{m-2}\Phi_{\cZ,\cZ_{\eta,\be}}(\hb,u,z),
\end{split}\end{equation*}
as claimed in~\e_ref{PhiZstr_e}.\\

\vspace{.2in}

\noindent
{\it Department of Mathematics, SUNY Stony Brook, NY 11794-3651\\
azinger@math.sunysb.edu}

\end{document}